\documentclass[11pt,twoside]{article}
\usepackage{latexsym}
\usepackage{amssymb,amsbsy,amsmath,amsfonts,amssymb,amscd}
\usepackage{epsfig, graphicx,subfigure}
\newcommand {\wt}[1] {{\widetilde #1}}

\newcommand{\commentout}[1]{}
\newcommand{\R}{\mathbb{R}}

\newcommand{\1}{{\mathchoice {\rm 1\mskip-4mu l} {\rm 1\mskip-4mu l}
{\rm 1\mskip-4.5mu l} {\rm 1\mskip-5mu l}}}
\newcommand{\ep}{\varepsilon}

\newcommand {\Da} {\Delta}

\newcommand {\lb} {\lambda}
\newcommand {\Chi} {{\bf \raise 2pt \hbox{$\chi$}} }

\newcommand {\cac} { {\mathcal C} }

\newcommand {\opG} { {\mathcal G} }
\newcommand {\f}   {\frac}
\newcommand {\p}   {\partial}
\newcommand {\grad}{\nabla}
\newcommand {\g}   {\bf}
\newcommand{\beq}{\begin{equation}}
\newcommand{\beqa} {\begin{array}{rl}}
\newcommand{\eeq}{\end{equation}}
\newcommand{\eeqa}{\end{array}}
\newtheorem{theorem}{Theorem}[section]
\newtheorem{lemma}[theorem]{Lemma}

\newtheorem{proposition}[theorem]{Proposition}

\newcommand{\qed}{{ \hfill
                       {\unskip\kern 6pt\penalty 500
                       \raise -2pt\hbox{\vrule\vbox to 6pt{\hrule width 6pt
                       \vfill\hrule}\vrule} \par}   }}
\title{\Large \bf Analysis of a Population Model Structured by the Cells Molecular Content}

\author{Marie Doumic \thanks{
D\'epartement de Math\'ematiques et Applications,
\'Ecole Normale Sup\'erieure, projet INRIA BANG, 
            45 rue d'Ulm, F~75230 Paris cedex 05, France;
email: doumic@dma.ens.fr} \thanks{INRIA Rocquencourt, projet BANG, Domaine de Voluceau, BP 105, F 78153 Rocquencourt, France; email: marie.doumic@inria.fr}
}

\date{\today}

\begin{document}
\maketitle
\pagestyle{plain}
\pagenumbering{arabic}

\begin{abstract}
We study the mathematical properties of a general model of cell division structured with several internal variables. We begin with a simpler and specific model with two variables, we solve the eigenvalue problem with strong or weak assumptions, and deduce from it the long-time convergence. The main difficulty comes from natural degeneracy of birth terms that we overcome with a regularization technique.  We then extend the results to the case with several parameters and recall the link between this simplified model and the one presented in \cite{CBBP1}; an application to the non-linear problem is also given, leading to robust subpolynomial growth of the total population. 
\end{abstract}

{\bf AMS subject classification:} 35A05 - 35P05 - 92D25 - 70K20
 
 \

{\bf Keywords:} structured populations, cell division, relative entropy, long-time asymptotic, eigenproblem, transport equation.

\

\section{Introduction}

\subsection{Presentation of the Model and Link with Other Models}\label{sec:pres}

Regulation of the cell division cycle governs the development of all organisms. Understanding it is central to the study of homeostasis, tumour growth and cancer,  but is made particularly difficult due to the numerous phenomena that can have an influence on it (see for instance \cite{Whitfield} and \cite{CMTU}, or \cite{CELL} for a general presentation of the cell cycle). 

For these reasons, several models have been proposed usually structured by a single variable (age, size, etc) \cite{MD}. But modern biology offers more accurate structuring variables as proteins or molecular content \cite{RCS}.

In order to investigate qualitatively the long-time behaviour of a cell population, 
we consider here a general model structured both in age, represented by the variable $a,$ and in another aggregated variable $x.$ Namely:
\begin{equation} \left\{
	\begin{array}{cl}
& \f{\p}{\p t}  n +  \f{\p}{\p a} n + \f{\p}{\p x}[ \Gamma(a,x) n ] + B(a,x) \, n =0, \quad \, a \geqslant 0, \; x \geqslant 0, 
\\ \\
&n(t,a=0,x) = 2 \int b(a,x, y) n(t,a,y) dy \, da,
\\ \\
&\Gamma_+(a,x=0)n(t,a,x=0)=0\quad \forall \; a\geqslant \; 0.
\end{array} \right.
\label{eq:cellcyclin:t1}
\end{equation}
This variable $x$ can represent one of the various proteins produced  (cf. for instance \cite{RYASE} or \cite{CBBP1}), the maturity of the cell (as in \cite{ACR} and \cite{APM} for instance), its size, its DNA content (as in \cite{BBMJBWW} for instance) etc.
Our study can be generalized to the case when $x \in \R_+^n,$ that is, when several phenomena influencing the cell cycle are taken into account (see part 
\ref{extension_multiple}) ; in part \ref{part:appli:2phase}, we also investigate the possible application of this model to a non linear two cell-compartment model (to model proliferating and quiescent cells). 

We have supposed here that age evolves like time, \emph{i.e.} $\f{da}{dt}=1.$ The function $\Gamma=\f{da}{dx}$ represents the rate at which the $x$ content of a cell increases with age.
The function $B$ represents the total division rate, and $b(a,x,y)$ is the repartition function of a mother cell of $y-$ content to a daugther cell of $x-$ content. We impose for consistancy (see \cite{CBBP1})
\beq
B(a,y)=\int b(a,x, y) dx, 
\label{as:b1}
\eeq
\beq y B(a,y)=2 \int x b(a,x, y) dx,\qquad b(a,x,y)=b(a,y-x,y),\qquad b(a,x>y,y)=0.\label{as:b2}
\eeq
It means conservation of the number of cells, conservation in the $x$ variable and symmetry of the division, when $x$ represents a molecular content. 

We have taken the coefficients independent of time (see \cite{CGP} for a model with time-periodic coefficients).
If we suppose that the coefficients $\Gamma$ and $B$ do not depend on $x$
and that for the solution $\Gamma n$ vanishes for $x=0$ and $x=\infty,$ integrating equation (\ref{eq:cellcyclin:t1}) in $x$ and denoting ${\mathcal N}(t,a)=\int n(t,a,x) dx$ gives:
\begin{equation} \left\{
	\begin{array}{cl}
& \f{\p}{\p t}  {\cal N} +  \f{\p}{\p a} {\cal N} + B(a) \, {\cal N} =0, \quad \, a \geqslant 0, 
\\ \\
&{\cal N}(t,a=0) = 2 \int B(a) {\cal N}(t,a) da.
\end{array} \right.
\label{eq:McK}
\end{equation}
This is the classical linear McKendrick-Von Foerster equation, with a death term which is exactly half the birth term (see \cite{BP1} for a complete study of this equation in the case whithout death term, and \cite{M:McKendrick} for extension of this equation to a non-linear case).

If we suppose that the coefficients $\Gamma,$ $B$ and $b$ do not depend on the age variable $a,$ that the integral ${\cal N}(t,x)=\int\limits_0^\infty n(t,a,x) da$ converges and that $\lim\limits_{a\to\infty} n(t,a,x) = 0,$ we can integrate (\ref{eq:cellcyclin:t1}) in $a$ and find:
\begin{equation} \left\{
	\begin{array}{cl}
& \f{\p}{\p t}  {\cal N} +  \f{\p}{\p x} [{\Gamma(x) \cal N}] + B(x) \, {\cal N}(x) =2 \int  b(x,y) {\cal N}(t,y) dy,
\\ \\
&\Gamma(x=0){\cal N}(t,x=0)=0\quad \forall \; t\geqslant \; 0.
\end{array} \right.\label{eq:size-structured}
\end{equation}
If $\Gamma = 1,$ we find the pure size-structured model, which has been studied in \cite{M:exist}, \cite{MMP}, \cite{PR},  \cite{PZ} for instance. In the case $\Gamma=1$ or $\Gamma=x^\mu,$ existence of a solution to the eigenvalue problem for a general $b$ is proved in \cite{M:exist} using approximation scheme. 

\

A model describing the dynamics of a cell population divided into proliferative and quiescent compartments was presented in \cite{CBBP1} and \cite{CBBP2}. It can be written:

\begin{equation} 
\label{eq:cellcyclin:t:2p}
\left\{
	\begin{array}{cl}
& \f{\p}{\p t} p + \f{\p}{\p a} p + \f{\p}{\p x}[ \Gamma(a,x) p ] + \big(B(a,x)+d_1+L(a,x)\big) \, p -  G(N(t)) q =0, \quad \, a \geqslant 0, \; x \geqslant 0, 
\\ \\
& \f{\p}{\p t} q + (G(N(t)) + d_2) q = L(a,x) p,\quad \, a \geqslant 0, \; x \geqslant 0, 
\\ \\
&p(t,a=0,x) = 2 \int b(a,x, y) p(t,a,y) dy \, da,
\\ \\
& p,\;q \geqslant 0, \qquad \int p(t,a,x)\,+\,q(t,a,x)\, dx \, da =N(t).
\end{array} \right.
\end{equation}
We will see in part \ref{part:appli:2phase} how it can be reduced to the study of problem (\ref{eq:cellcyclin:t1}).

\subsection{The Eigenvalue Problem}
\label{part:eigenproblem}
In order to study the asymptotic behaviour of 
the solution of problem (\ref{eq:cellcyclin:t1}), we consider the eigenvalue problem (see \cite{CBBP1} and \cite{CBBP2}): find $(\lb_0, N)$ solution to
\label{sec:model}

%
\medskip
\begin{equation} \left\{
	\begin{array}{cl}
& \f{\p}{\p a} N + \f{\p}{\p x}[ \Gamma(a,x) N ] + \big(\lb_0+B(a,x)\big) \, N =0, \quad \, a \geqslant 0, \; x \geqslant 0, 
\\ \\
&N(a=0,x) = 2 \int b(a,x, y) N(a,y) dy \, da,
\\ \\
& N \geqslant 0, \qquad \int N dx \, da =1.
\end{array} \right.
\label{eq:cellcyclin}
\end{equation}
This is an original problem which can be seen as a usual Cauchy problem where the initial data is related to the ``future''.

It is useful also to study the adjoint problem: 
\begin{equation} \left\{
	\begin{array}{cl}
& -\f{\p}{\p a} \phi  -   \Gamma(a,x)  \f{\p}{\p x }\phi + \big(\lb_0+B(a,x)\big) \, \phi =2 \int  b(a,y,x) \phi(0,y) dy, \quad \, a\geqslant 0, \; x \geqslant 0, 
\\ \\& \phi \geqslant 0, \qquad \int \phi N dx \, da =1.
\end{array} \right.
\label{eq:cellcyclind}
\end{equation}

Here we make the assumptions:

\beq\label{as:gamma0}
\Gamma(a,x)= x \wt \Gamma (a,x),
\eeq 
with 
\beq \left\{
\beqa
\wt \Gamma (a,x)& > 0 \quad \text{for }  x\approx 0,
\\ \\
\wt \Gamma (a,x) &\leqslant 0 \quad \text{for }  x \geqslant  x_M.
\eeqa 
\right. 
\label{as:gamma1}
\eeq
Then there is no need of boundary condition at $x=0$ and conservation in the $x$ variable is enforced according to the biophysical interpretation, when $x$ represents a molecular content. The fact that $\Gamma$ becomes negative beyond a maximum value $x_M$ means that the $x-$content of the cells remains bounded, but the results can be generalised to the case when $\Gamma$ remains nonnegative everywhere.

As a consequence of condition (\ref{as:b1}), integrating equation (\ref{eq:cellcyclin}), 
 we have for $N$ vanishing at infinity
$$
\lb_0 = \int B(a,x) N(a,x) dx \, da,
$$
(in words, the population number can only grow by cell division).
Integrating the equation (\ref{eq:cellcyclin}) against the weight $x,$ as soon as
$\lim\limits_{a\to\infty} \int xN(a,x) dx = 0,$
we have, using (\ref{as:b2})
$$
\lb_0 \int  x N(a,x) dx \, da= \int \Gamma(a,x) N(a,x) dx \, da ,
$$
(in words, the total molecular content can only increase by the reaction term $\Gamma$).
Integrating the equation against the weigth $a,$ as soon as $\lim\limits_{a\to\infty} \int a N(a,x) dx = 0,$
we have
\beq
\lb_0 \int  a N(a,x) dx \, da+ \int a B(a,x) N(a,x) dx \, da=1.
\eeq
Also, we can do as in \cite{M:exist} and for $0<\eta<1,$ if 
$\lim\limits_{a\to\infty} \int \, e^{\lb_0 \eta a} N(a,x) dx =0,$
multiplying by $e^{\lb_0 \eta a}$ and integrating, we find
$$\forall \; 0<\eta < 1, \quad \int e^{\lb_0 a \eta} N(a,x) dx da \leqslant \f{1}{1-\eta}.$$


\

We can reduce the study to the solutions on the domain $(a,x)\in \R_+\times [0,x_M].$  Indeed, 
using the method of characteristics based on the solution to the differential system parametrized by the Cauchy data $(a,x):$
\beq \left\{
\begin{array}{cl} &
\f{d}{d a} X(a,x)= \Gamma \big(a, X(a,x)\big),\quad a\geqslant 0,\; x\geqslant 0, 
\\ \\
& X(0,x)=x,\quad x\geqslant 0,
\end{array}
\right. 
\label{eq:char}
\eeq
Cauchy-Lipschitz theorem gives us, as soon as 
$\Gamma \in \cac^1_b (\R_+\times\R_+)$
for instance, the existence and uniqueness of the flow $X(a,x).$ We denote $Y(a,x)$ the inverse flow, defined by 
\beq
Y\big(a,X(a,y)\big)=y.
\label{eq:chard}
\eeq

We deduce from (\ref{as:gamma1}) that for all
$x\leqslant x_M,$ for all $a\geqslant 0,$ $X(a,x) \leqslant X (a,x_M) \leqslant x_M.$
The formula (\ref{eq:char:N}), proved below in lemma \ref{lem:reform}, shows that if the solution $N$ verifies  $N(a=0,\; x\geqslant x_M)=0,$ then $N\big(a\geqslant 0,\; x\geqslant X(a,x_M)\big)=0,$ and in particular $N(a,x \geqslant x_M)=0.$ Thus,
we add the following condition to problem (\ref{eq:cellcyclin}):
\beq\label{as:suppN}
N(0,x\geqslant x_M)=0,
\eeq
and it allows us to restrict our study to the compact set $[0,\,x_M].$
We could also exchange assumption (\ref{as:gamma0}) with the following one:
\beq 
\label{as:N0}
N(a\geqslant 0,x=0)=0,\qquad b(a\geqslant 0,0,y)=0.
\eeq
Contrarily to the solution $N,$ the solution of the adjoint problem (\ref{eq:cellcyclind}) does not necessarily have its support in $[0,x_M],$ but we make the following assumption on the support of the function $b$ (cf. the proof in section \ref{sec:proofs}):
\beq\label{as:suppb}
b\big(a,x,y\geqslant X(a,x_M)\big)=0,
\eeq
which implies that 
$b(a,x>x_M,y)=0,$ and thus, by formula (\ref{form:N}) proved below, it implies (\ref{as:suppN}).



In all the following, if nothing is specified, we suppose assumptions (\ref{as:gamma0}), (\ref{as:gamma1}), (\ref{as:b1}) and (\ref{as:b2}) are satisfied. We denote $X(a,x)$ the characteristic flow solution of (\ref{eq:char}) and $Y(a,x)$ the inverse flow defined by (\ref{eq:chard}). 

\subsection{Reformulation of the Problem with the Method of Characteristics}

We first give the following formulae, on which the proofs are based.

\begin{lemma}\label{lem:reform}
For $\Gamma\in \cac^1 (\R_+ \times [0,x_M]),$ under assumption (\ref{as:suppN}) a solution $N$ to (\ref{eq:cellcyclin}) verifies the following formula (as soon as the integral converges):
\beq\label{form:N}
 N (0,x) = 2 \int\limits_0^\infty \int \limits_0^{x_M} b\big(a,x,X(a,y)\big)
 N (0, y)e^{ -\int\limits_0^a \big\{\lb_0 + B \big(s,X(s,y)\big) \big\}ds}\, dy da,
\eeq
and we can write:
\beq\label{eq:char:N1}
N\big(a,X(a,x)\big)=N(0,x)e^{ -\int\limits_0^a \big\{\lb_0 + \f{\p}{\p x} \Gamma \big(s,X(s,x)\big) + B \big(s,X(s,x)\big) \big\}ds},\eeq
or also
\beq\label{eq:char:N}
N(a,y)=N(0,Y(a,y))e^{ -\int\limits_0^a \big\{\lb_0 + \f{\p}{\p x} \Gamma \big(s,X\big(s,Y(a,y)\big)\big) + B \big(s,X\big(s,Y(a,y)\big)\big) \big\}ds}.\eeq
Under assumption (\ref{as:suppb}) a solution $\phi$ to (\ref{eq:cellcyclind}) verifies (as soon as the integral converges): 
\beq\label{form:phi}
 \phi (0,x) = 2\int\limits_0^\infty\int\limits_0^{x_M}  b\big(a,y,X(a,x)\big) \phi (0,y) e^{ -\int\limits_0^a \big\{\lb_0 +  B \big(s,X(s,x)\big) \big\}ds}dy da,
\eeq
and as soon as the integral converges, we can write:
\beq\label{eq:phi}
\phi(a,x)= 2  \int\limits_a^\infty \int\limits_0^{x_M}  b\big(s,y,X\big(s,Y(a,x)\big)\big) \phi(0,y) e^{-\lb_0(s-a)-\int\limits_a^s B\big(\sigma,X(\sigma,Y(a,x))\big)d\sigma} dy \, ds.
\eeq
\end{lemma}

\begin{proof}
We set
$\wt N  (a,x)= N \big(a, X(a,x)\big)e^{ \int\limits_0^a \big\{\lb_0+ \f{\p}{\p x} \Gamma \big(s,X(s,x)\big) + B \big(s,X(s,x)\big) \big\}ds}$
and rewrite equation (\ref{eq:cellcyclin}):
$$
\f{\p}{\p a} \wt N (a,x) = \big[ \f{\p}{\p a} N + \f{\p}{\p x}[ \Gamma N ] + \big(\lb_0+B \big) \, N \big] \big(a, X(a,x)\big)
e^{ \int\limits_0^a \big\{\lb_0 + \f{\p}{\p x} \Gamma \big(s,X(s,x)\big) + B \big(s,X(s,x) \big\}ds\big)}
= 0.$$
This gives the equalitities (\ref{eq:char:N1}) and (\ref{eq:char:N}). 
We can rewrite the boundary condition of problem (\ref{eq:cellcyclin}) as: 
$$ N (0,x) = 2 \int\limits_0^\infty \int\limits_0^{x_M}
b (a,x,y)
 N \big(0, Y(a,y)\big)e^{ -\int\limits_0^a \big\{\lb_0 + \f{\p}{\p x} \Gamma \big(s,X\big(s,Y(a,y)\big)\big) + B \big(s,X\big(s,Y(a,y)\big)\big) \big\}ds}dy \, da,$$
or equivalently, changing variables $y \to Y(a,y)$ and noting that
$
e^{-\int\limits_0^a \f{\p}{\p x} \Gamma \big(s,X(s,y)\big) 
ds}= \f{1}{\p_y X(a,y)},$
we find formula (\ref{form:N}).

In the same way, we set
 $$ \wt \phi  (a,x)= \phi \big(a, X(a,x)\big)e^{- \int\limits_0^a \big\{\lb_0  + B \big(s,X(s,x)\big) \big\}ds}$$
and rewrite equation (\ref{eq:cellcyclind}) as
$$\begin{array}{cl} \f{\p}{\p a} \wt \phi (a,x) & = \big[ \f{\p}{\p a} \phi + \Gamma \f{\p}{\p x} \phi - \big(\lb_0+B \big) \, \phi \big] \big(a, X(a,x)\big)
e^{ -\big(\int\limits_0^a \big\{\lb_0  + B \big(s,X(s,x)\big) \big\}ds\big)}
\\
& = -2 \int\limits_0^{x_M}  b\big(a,y,X(a,x)\big) \phi(0,y)  e^{ -\int\limits_0^a \big\{\lb_0 + B \big(s,X(s,x)\big) \big\}ds}\,dy.
\end{array}
$$
Integrating along $a$ on $\R_+$ we get formula (\ref{form:phi}), and integrating from $0$ to $a$ we then find formula (\ref{eq:phi}).
\end{proof}

\subsection{Some Examples}\label{sec:examples}

\

{\bf Case 1:} for the function $\Gamma,$ we can take for instance
$\Gamma (a,x) = C_1 x (x_M - x).$
We easily calculate that
$$X(a,x) =\f{x\cdot x_M e^{x_M C_1 a}}{x_M + x(e^{x_M C_1 a} -1)},$$
so for $x \approx 0,$ one has $X(a,x) \approx x e^{x_M C_1 a}.$

{\bf Case 2:} another possible example is to take
$\Gamma (a,x) = C_1 x^\alpha (x_M - x)^\beta,$ with $\quad 0<\alpha <1$ and $\beta >0.$
For $x>0,$ we then have $$X(a,x) \approx_{x\approx 0} (Ca+ x^{1-\alpha})^{\f{1}{1-\alpha}}.$$

\begin{figure}[ht]
\begin{center}
\includegraphics[width=7cm]{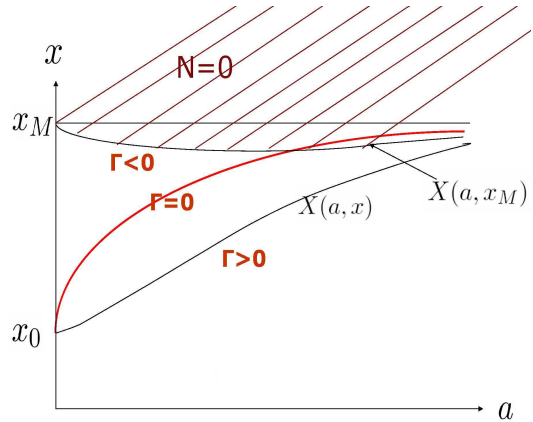}
\caption{\label{fig:gamma} Form of $\Gamma$ defined by (\ref{def:gamma}), and its related characteristics.}
\end{center}
\end{figure}

{\bf Case 3:} biological considerations lead \cite{BBCA} and \cite{CBBP1} to take for $\Gamma$ (case illustrated in fig. \ref{fig:gamma}):
\beq \label{def:gamma}
\Gamma(a,x)=c_1 \f{x}{1+x} \bigl(r_1-r_2 e^{-c_4a}\bigr)-c_2 x,\qquad r_1,\;r_2 >0,\qquad\f{c_2}{c_1} <r_1-r_2.
\eeq 
In this case, $x_M= \f{c_1}{c_2} r_1 -1$ and $\Gamma(0,x_0)=0$ with $x_0=\f{c_1}{c_2} (r_1 - r_2)-1 >0.$ 

For the division rate $B,$ we can take
\beq\label{eq:defB}
B(a,x)=C_2 x^\gamma\1_{A^*\leqslant a\leqslant A_1},\quad \gamma \geqslant 1,\quad 0\leqslant A^* <A_1\leqslant \infty.
\eeq
To define the repartition function $b$ (see \cite{HLFP}, \cite{SMPT} and \cite{SY} for biological motivations) we can choose a uniform repartition:
$$b(a,x,y)=\1_{x\leqslant y} \f{B(a,y)}{y}.$$
On the opposite, the classical example of equal mitosis is given by $$b(a,x,y)= \delta_{x=\f{y}{2}} B(a,y).$$
It is the case for instance if $x$ represents the DNA-content of the cells.
In this last case, formulae (\ref{form:N}) and (\ref{form:phi}) of lemma \ref{lem:reform} are no longer valid; see part \ref{dirac} for an adaptation of the proof. A study of the link between the repartition function $b$ and the proliferation rate $\lb_0,$ for the pure size-structured model, can be found in \cite{M:optim}.

\noindent
We can also generalize the uniform repartition by
$b(a,x,y)=\f{B(a,y)}{y(1-2\eta)}\1_{\eta y\leqslant x\leqslant (1-\eta) y}$ with $0<\eta<\f{1}{2},$ or take a Poisson distribution for $b.$

\section{Resolution of the Problem in a Regular Case}

\label{part:strong}

			\subsection{Main Results}\label{sec:main}
In this part, we consider regular functions $b$ and $B$ in order to obtain strong compactness and regular solutions with non extinction, \emph{i.e.}, $\lb_0 \geqslant 0.$ We need the following assumptions.



\beq\label{as:gamma3} 
\Gamma \in \cac^1_a (\R_+, \cac^2_x (\R_+)).
\eeq
This condition is used to prove the regularity of the solution $N$ (else it can be weakened: see part \ref{part:gen:weak}). 

\beq\label{as:B1}
B \in \cac \big(\R_+\times [0,x_M]\big),\quad \f{\p}{\p x} B\in \cac \big (\R_+ \times [0,x_M]\big),\quad B(a,0)=0,
\eeq




\beq\label{as:b4}
b  \in L^\infty (\R_+^3), \quad b(a,.,y) \in \cac^1_b\big(\{x \in [0,x_M]\; ; \quad x\leqslant y\big\}\big).
\eeq
The function $b$ is not necessarily continuous for $x\in\R_+$ since $b(a,x>y,y)=0,$ but it is necessary to suppose it is regular for $x\leqslant y$ in order to prove strong compactness. 
On the previous examples, except the case of equal mitosis, and with $\gamma\geqslant 1,$ the function $b$ verifies assumption (\ref{as:b4}).



\beq\label{as:H1}
\int\limits_0^\infty \int\limits_0^{x_M} e^{-\int\limits_0^a B\big(s,X(s,x)\big)ds} dx\, da <\infty.
\eeq 
This is a key assumption, which is used to obtain compactness and a solution $N$ in $L^1(\R_+^2).$ We can give an equivalent formulation in terms of partial differential equations, as follows.

\noindent
Let $v$ the solution of the following problem:
\beq\label{as:H1:PDE}
\begin{array}{ll}
&\f{\p}{\p a} v + \f{\p}{\p x}[ \Gamma(a,x) v ] + B(a,x) \, v =0, \quad \, a \geqslant 0, \;x_M\geqslant  x \geqslant 0, 
\\ \\ 
&v(a=0,x) = 1\quad 0\leqslant x\leqslant x_M,
\end{array}
\eeq
then (\ref{as:H1}) means $v\;\in L^1(\R_+^2).$ 
Assumption (\ref{as:H1}) can be replaced by slightly more general ones, allowing in particular $\int B\bigl(a,X(a,x)\bigr)da$ to be finite: see for instance paragraph \ref{part:gen:weak}.

Finally, we need the following two assumptions to prove uniqueness of a solution.
\beq 
\label{as:Bpositiv}
\exists A_1\geqslant A^*> 0,\;\forall\;A^*\leqslant a\leqslant A_1,\;\forall\;y\;\in\;]0,x_M[,\quad B(a,y)>0,
\eeq 

\beq
\label{as:bpositiv}
\forall \; x\; \in \; ]0,\; x_M[,\quad \forall\; y\;\in\;]0,x_M[,\quad \int b\big(a,x,X(a,y)\big) da \; >0.
\eeq
Assumption (\ref{as:bpositiv}) implies that $Supp_a(b)$  is unbounded and that
$\lim\limits_{a\to\infty} X(a,0<y<x_M)=x_M.$

\noindent
Uniqueness could also be proved under other assumptions than (\ref{as:Bpositiv}) and (\ref{as:bpositiv}): see for instance part \ref{part:compact}.

\

\noindent
In the previous examples, the function $\Gamma$ verifies assumption (\ref{as:gamma3}) but we need to regularize  $B$ in order to obtain   (\ref{as:B1}), and also to regularize $b,$ in the case of equal mitosis, to ensure (\ref{as:b4}). The positivity conditions (\ref{as:Bpositiv}) and (\ref{as:bpositiv}) stand if $A_1=+\infty$ in the definition of $B,$ but assumption (\ref{as:bpositiv}) is not true in the case of equal mitosis.  It remains to check assumption (\ref{as:H1}).
If $\Gamma(a,x)=C_1 x (x_M-x),$ we have to prove the convergence of
$\iint e^{-\int\limits_0^a C_2 X(s,x)^\gamma ds} dxda,$
or equivalently of 
$$\iint \int e^{-x^\gamma e^a} dxda 
=
\iint \f{e^{-x^\gamma u}}{u} dxdu
=
\iint \f{e^{-X^\gamma}}{u^{1+\f{1}{\gamma}}} dXdu.
$$
Since $\gamma>0,$ this integral converges, which gives us assumption (\ref{as:H1}).

If $\Gamma (a,x) \approx_{x\approx 0} C x^\alpha$ with $0<\alpha <1,$ 
we have to study the convergence of
$$I=
\int\limits^\infty\int\limits_0 e^{-\int\limits_0^a (Cs+x^{1-\alpha})^{\f{\gamma}{1-\alpha}} ds} dx\, da 
\approx C \int\limits^\infty e^{-a^{1+ \f{\gamma}{1-\alpha}}} da.$$
Since $1+ \f{\gamma}{1-\alpha}>0,$ this integral converges.



\begin{theorem}\label{th:cellcyclin}
Under assumptions (\ref{as:suppN}), (\ref{as:gamma3})--(\ref{as:H1}), there exists a solution $N \in \cac^1_b\big( \R_+^2\big)$ to (\ref{eq:cellcyclin}) for a $\lb_0 > 0.$ Moreover, $N\big(a,x\geqslant X(a,x_M)\big)=0.$
\end{theorem}



\begin{theorem} \label{th:cellcyclind}
Under the assumptions of theorem \ref{th:cellcyclin}, with the additional assumption (\ref{as:suppb}), 
 there exists a solution $\phi \in \cac^1_b (\R_+^2)$ of (\ref{eq:cellcyclind}) 
for a $\lb_0 > 0.$ 
\end{theorem}



\begin{theorem} \label{th:uniqueness}
Under the assumptions of theorem \ref{th:cellcyclin} and with the additional assumptions (\ref{as:suppb}), (\ref{as:Bpositiv}) and (\ref{as:bpositiv}), the solutions $(\lb_0,\; N,\;\phi)$ of the eigenproblems (\ref{eq:cellcyclin}) and (\ref{eq:cellcyclind}) are unique.
\end{theorem}


					\subsection{Regularised Problem and Method of Characteristics}

\label{part:char}

In order to apply Krein-Rutman theorem (refer to \cite{DL} for instance), we first have to consider a regularised problem, where the operator is strictly positive. 
  We write
\beq \label{def:bBep:strong}
b_\ep (a,x,y) = b (a,x,y) +  \f{\ep}{x_M},\quad B_\ep (a,y)= B(a,y) + \ep,
\eeq
for $\ep >0.$
The functions $b_\ep$ and $B_\ep$ verify the same relation (\ref{as:b1}) than $b$ and $B$. We  consider the regularised problem:
\begin{equation} \left\{
	\begin{array}{cl}
& \f{\p}{\p a} N_\ep + \f{\p}{\p x}[ \Gamma(a,x) N_\ep ] + \big(\lb+B_\ep(a,x)\big) \, N_\ep =0, \quad \, a \geqslant 0, \; 0\leqslant x \leqslant x_M, 
\\ \\
&N_\ep(a=0,x) = 2 \int b_\ep(a,x, y) N_\ep(a,y) dy \, da, \quad 0\leqslant x\leqslant x_M,
\\ \\
& N_\ep \geqslant 0, \qquad \int N_\ep dx \, da =1,
\end{array} \right.
\label{eq:cellcyclin:ep}
\end{equation}
 and its adjoint:
\begin{equation} \left\{
	\begin{array}{cl}
& -\f{\p}{\p a} \phi_\ep  -   \Gamma(a,x)  \f{\p}{\p x }\phi_\ep + \big(\lb+B_\ep(a,x)\big) \, \phi_\ep =2 \int \phi_\ep(0,y) b_\ep(a,y,x) dy, \quad \, 0\leqslant a, \; 0\leqslant x \leqslant x_M, 
\\ \\
& \phi_\ep \geqslant 0, \qquad \int \phi_\ep N_\ep dx \, da =1.
\end{array} \right.
\label{eq:cellcyclind:ep}
\end{equation}
We can write the formula (\ref{form:N}) for this regularised problem:
\beq
 N_\ep (0,x) = 2 \int\limits_0^\infty \int \limits_0^{x_M} b_\ep \big(a,x,X(a,y)\big)
 N_\ep (0, y)e^{ \big(-\int\limits_0^a \big\{\lb + B_\ep \big(s,X(s,y)\big) \big\}ds\big)}dy \, da.
\eeq
To find a solution to this equation, we consider the Banach space
$$E= \cac^1_b \big([0, x_M], \R\big),\quad ||f||_E:=||f||_{L^\infty} + ||\f{d}{dx} f||_{L^\infty},$$ 
and study the linear integral operator $\opG_{\lb}^\ep$ defined on $E$ by:
\beq 
\label{def:opG}
\opG_{\lb}^\ep(f)(x)= 2 \int b_\ep \big(a,x,X(a,y)\big)
 f(y)e^{ -\int\limits_0^a \big\{\lb + B_\ep \big(s,X(s,y)\big) \big\}ds}dy \, da.
\eeq


\begin{lemma}\label{lem:comp}
Under the assumptions of theorem \ref{th:cellcyclin},  
for all $\ep \geqslant 0,$ $\lb \geqslant 0,$ the operator $\opG^{\ep}_\lb$ is compact.
\end{lemma}


\begin{proof}
In order to apply the Ascoli theorem, we evaluate, for an arbitrary $f\in E,$ $||f||_E\leqslant 1,$ the quantities  
$\Da(x_1,x_2)=\opG_{\lb}^\ep(f)(x_1)-\opG_{\lb}^\ep(f)(x_2)$ and  $\Da'(x_1,x_2)=\f{d}{dx}\opG_{\lb}^\ep(f)(x_1)-\f{d}{dx}\opG_{\lb}^\ep(f)(x_2),$ 
and prove they tend to zero uniformally on $x_1, \, x_2$ when $|x_1 - x_2|$ tends to zero. 
Supposing $x_1 < x_2,$ we can write:
$$\Da\leqslant 2 \int\limits_0^\infty e^{-\int\limits_0^a B\big(s, X(s,y)\big)ds }\biggl\{
\int\limits_{Y(a,x_2)}^{x_M}| b \big(a,x_1,X(a,y)\big)-b\big(a,x_2,X(a,y)\big)| 
+ \int\limits_{Y(a,x_1)}^{Y(a,x_2)}b\big(a,x_1,X(a,y)\bigl)\biggr\} dyda
.$$

For each term, thanks to assumption (\ref{as:H1}) and to the fact that $b\in L^\infty (\R_+^3),$ we first restrict the integral to $a\leqslant A$ where $A$ is chosen such that  the rest $\int\limits_A^\infty$ be sufficiently small independantly of $x_1$ and $x_2$. The first term may be bounded  by $C|x_1-x_2|$ thanks to (\ref{as:b4}). The second term may be bounded as $O(|x_1 - x_2|)$ since $Y(a,x)$ is equicontinuous on $[0,A]\times[0,x_M].$
We evaluate $\Da'(x_1,x_2)$ in the same way, with the help of $||\f{d}{dx}f||_{L^\infty\big([0,x_M]\big)}$. 
This proves the equicontinuity of the family $\opG^\ep_\lb (f),$ for $f\in E$ and $||f||_E\leqslant 1,$ therefore $\opG^\ep_\lb$ is compact.

\end{proof}
\



\begin{lemma}\label{lem:reg}
For all $\ep >0$ and $\lb \geqslant 0,$ under the assumptions of theorem \ref{th:cellcyclin}, there exists a unique $\mu_{\lb,\ep} >0$ and a unique $N^0_{\lb,\ep}\in E,$ $N^0_{\lb,\ep} >0,$ such that 
		$$\opG_{\lb}^\ep(N^0_{\lb,\ep})=\mu_{\lb,\ep} N^0_{\lb,\ep},\quad
||N^0_{\lb,\ep}||_E = 1.
$$
Moreover, $\lb \to \mu_{\lb,\ep}$ is a continuous decreasing function, vanishing  when $\lb$ tends to infinity and taking the value $2$ when $\lb=0.$
\end{lemma}



\begin{proof}
The compact operator $\opG_{\lb}^\ep$  is strictly positive on E, thanks to 
the fact that $\ep >0$,
so Krein-Rutman theorem gives existence and uniqueness of $\mu_{\lb, \ep}>0$ and $N^0_{\lb,\ep}>0.$
The function $\lb \to \mu_{\lb,\ep}$ decreases because $\opG_\lb^{\ep}$ decreases with $\lb.$
It is continuous thanks to the uniqueness of the eigenvector and to the compactness of the family of operators $\opG_\lb^\ep.$
For $f\in E,$ $f>0$ and $||f||_E = 1,$ we have:
$$\opG_{\lb}^\ep(f)(x)  \leqslant  2 \int  f(y) e^{-\lb a} ||b_\ep||_\infty dy\, da = \f{2 ||b||_\infty +2\ep/x_M }{ \lb} \int f(y) dy.$$
We evaluate this quantity in $f=N^0_{\lb,\ep},$ and integrate in $x:$
\beq\label{ineq:lb}
\mu_{\lb,\ep} \leqslant \f{2||b||_\infty x_M +2\ep}{\lb},
\eeq
so $\lim\limits_{\lb\to\infty} \mu_{\lb,\ep} = 0.$ 
Integrating in $x,$ we find $$\int \opG_{\lb}^\ep(f)(x) dx = 2 \int f(y) e^{-\lb a} \f{d}{da}\biggl(-e^{-\int\limits_0^a B_\ep\big(s,X(s,y)\big) ds}\biggr) dy\, da.$$
Integrating by parts, we get: $$\int \opG_{ \lb}^\ep(f)(x) dx = 2\int f(y) dy  -2 \int f(y) \lb e^{-(\lb+\ep) a} e^{-\int\limits_0^a B \big(s,X(s,y)\big) ds} dy\, da,
$$
so for $\lb=0,$ $\int \opG_{ \lb}^\ep(f)(x) dx = 2\int f(x) dx ,$  so $\mu_{0,\ep}=2.$

\end{proof}


\subsection{Proofs of the Main Theorems}
\label{sec:proofs}
Thanks to lemma \ref{lem:reg}, for all $\ep>0,$ we can define $\lb_\ep>0$ as being the only $\lb >0$ such that $\mu_{\lb_\ep,\ep}=1.$ We denote $N^0_\ep$ the associated eigenvector, with $||N^0_\ep ||_E =1.$ Since
$N^0_\ep =\opG^\ep_{\lb_\ep} (N^0_\ep),$
the family $(N^0_\ep)_{0\leqslant\ep \leqslant 1}$
is compact in $E.$ Indeed, as in lemma \ref{lem:comp}, we can show that it is an equicontinuous family of functions and apply the Ascoli theorem.
Thanks to inequality (\ref{ineq:lb}), we know that $\lb_\ep$ is bounded so
we can extract a subsequence $(\lb_\ep, N^0_\ep)$ tending to $(\lb_0,N^0)\in \R_+ \times E,$ $N^0\geqslant 0,$ $\lb_0 \geqslant 0$ and $||N^0||_E = 1.$ 
%
We have
$$N^0_\ep(x)= 2 \int \big\{b \big(a,x,X(a,y)\big)
                     N^0_\ep(y)
+\,\f{\ep}{x_M} 
         N^0_\ep(y)\big\}e^{ -\int\limits_0^a \big\{\lb_\ep +\ep + B \big(s,X(s,y)\big) \big\}ds}dyda.
$$
Under assumption (\ref{as:H1}), the second term in this expression vanishes when $\ep\to 0,$ so we have: 
\beq\label{eq:opG}
N^0(x)=\opG^0_{\lb_0} (N^0) (x)= 2 \int b \big(a,x,X(a,y)\big)
 N^0(y)e^{ -\int\limits_0^a \big\{\lb_0 + B \big(s,X(s,y)\big) \big\}ds}dy \, da.
\eeq
We define $N$ by
\beq\label{eq:N}
N(a,y)= \alpha N^0 ( Y(a,y)\big)e^{ -\int\limits_0^a \big\{\lb_0 + \f{\p}{\p x} \Gamma \big(s,X\big(s,Y(a,y)\big)\big) + B \big(s,X\big(s,Y(a,y)\big)\big) \big\}ds}\1_{y\leqslant X(a,x_M)}.
\eeq
To state that the so-defined function $N$ satisfies (\ref{eq:cellcyclin}), it remains only to prove that $\int N da dx <\infty,$ then we choose $\alpha$ such as $\int N da dx =1.$ We have
$\int N(a,y) da dy= \alpha \int N^0 (y)e^{ -\int\limits_0^a \big\{\lb_0 + B \big(s,X(s,y)\big) \big\}ds}dyda.$

\noindent
Thanks to assumption (\ref{as:H1}) and since $N^0\in L^\infty (\R_+),$ this integral converges. As seen in part \ref{part:eigenproblem}, integration of equation (\ref{eq:cellcyclin}) gives
$\lb_0=\int B N dx\, da >0,$
because  if it were zero, that would imply $\int\limits_0^{x_M} N^0(x) dx = 2 \int BN dx\, da =0,$ so $N^0=0,$ which is absurd since  $||N^0||_E=1.$
The continuity and the continuous derivability of $N$ at the points $y=X(a,x_M)$ come straigthforward from formulae (\ref{eq:opG}) and (\ref{eq:N}), and from the assumptions of regularity (\ref{as:gamma3})--(\ref{as:b4}): it ends the proof of theorem \ref{th:cellcyclin}.

\

In the same way than for the direct problem, formula (\ref{form:phi}) leads us to study the integral operator defined on the Banach space $E$ by 
$$\opG^{0\, *}_{\lb_0} (f) (x)=2\int\limits_0^\infty\int\limits_0^{x_M} b \big(a,y,X(a,x)\big) f (y) e^{ -\big(\int\limits_0^a \big\{\lb_0 +  B \big(s,X(s,x)\big) \big\}ds\big)} da dy.$$

The previous study for $\opG^0_{\lb_0}$ can be carried out in the same way, and we find an eigenvalue $\lb_1>0$ and 
an associated eigenvector $\phi^0 \in E$ for the operator $\opG^{0\,*}_{\lb_0}.$ We choose $||\phi^0||_{E} =1,$ and define $\phi(a,x)$ by the formula (\ref{eq:phi}).
It satisfies (\ref{eq:cellcyclind}), and it only remains to check that $\int \phi N da dx <\infty.$ We have:
$$\int \phi N da dx = 4 \alpha \int\limits_0^\infty \int\limits_0^{x_M}\int\limits_a^\infty \int\limits_0^{x_M} \phi^0(y) b\big(s,y,X\big(s,Y(a,x)\big)\big)e^{-\lb_1 a -\lb_0( s-a) }$$
$$ N^0\big(Y(a,x)\big) e^{- \int\limits_0^s B\big(\sigma, X\big(\sigma,Y(a,x)\big)\big) d\sigma-\int\limits_0^a \p_x \Gamma \big(\sigma, X\big(\sigma, Y(a,x)\big)\big)d\sigma} \,dx\,ds\,dy\,da.
$$
Changing variables $x\rightarrow Y(a,x)$ as we did in part \ref{part:char} and integrating according to the $s-$variable before the $a-$variable, one gets:
$$\int \phi N da dx = 4 \alpha \int\limits_0^\infty \int\limits_0^{x_M} \int\limits_0^{x_M} \phi^0(y) N^0 (x)b\big(s,y,X(s,x)\big)(\f{e^{-\lb_0 s} - e^{-\lb_1}s}{\lb_1-\lb_0})e^{- \int\limits_0^s B\big(\sigma, X(\sigma,x)\big) d\sigma}\, dx \, dy \, ds.$$
This integral converges since $\lb_0,\lb_1>0$ (we do not know yet if $\lb_0=\lb_1$; if it is the case we just replace $\f{e^{-\lb_0 s} - e^{-\lb_1}s}{\lb_1-\lb_0}$ by $se^{-\lb_0 s}$).
We choose $\alpha$ such that $\int \phi N dx\, da =1.$ It ends the proof of theorem \ref{th:cellcyclind}.

\

For uniqueness, assumption (\ref{as:bpositiv}) implies that for all
$x\;\in\;]0,x_M[,$ we have $N(0,x)>0$ and $\phi(0,x)>0.$
Thanks to Fredholm alternative, uniqueness of a solution to problem (\ref{eq:cellcyclin}) implies that to problem (\ref{eq:cellcyclind}), so we only have to prove uniqueness for one of that two problems.
Consider two solutions $(\lb_1,N_1)$ and $(\lb_0,\phi)$ of (\ref{eq:cellcyclin}) and (\ref{eq:cellcyclind}).  
Then we have: 
$$\begin{array}{cl}
&-\lb_1 \int N_1 \phi \,dx\, da = 
 \int \biggl(\f{\p}{\p a} N_1 + \f{\p}{\p x}[ \Gamma(a,x) N_1 ] + B(a,x)\, N_1 \biggr) \,\phi \, dx\, da
\\
& = \int \biggl(-\f{\p}{\p a} \phi  -   \Gamma(a,x)  \f{\p}{\p x }\phi + B(a,x) \, \phi - 2 \int \phi(0,y) b(a,y,x) dy\biggr) N_1 dx\, da 
=- \lb_0 \int N_1 \phi\, dx\, da.
\end{array} 
$$
Thanks to assumption (\ref{as:bpositiv}), we have
$\int N_1 \phi \,dx \,da >0,$ and it implies $\lb_1=\lb_0.$ 

\noindent
Let us now suppose two solutions $N_0$ and $N_1$ of problem (\ref{eq:cellcyclin}) for the same eigenvalue $\lb_0.$ We use the following lemma (for a proof, we refer to \cite{BP1}, proposition 6.3.)
 \begin{lemma}\label{lem:absN}
 Let $f\in {\cal{S}} ' (\R_+^2)$ be a solution of
 \beq\label{eq:cellcyclin:1}
 \f{\p}{\p a} f + \f{\p}{\p x}[ \Gamma(a,x) f ] + \big(\lb_0+B(a,x)\big) \, f =0, \quad \, a \geqslant 0, \; x \geqslant 0, \eeq
 then $|f|$ satisfies the same equation. 
 \end{lemma}
We set $\wt N = |N-N_1.|$ According to lemma \ref{lem:absN}, $\wt N$ verifies equation (\ref{eq:cellcyclin:1}). We take $\phi$ as a test function:
$$2\iiint \phi\,(0,x) b(a,x,y) |N-N_1| (a,y)  dx \,dy\,da = 2 \int \phi(0,x) \big| \iint b(a,x,y) (N-N_1)(a,y) dx\, dy\,\big| da. 
$$
Thanks to the fact that $\phi(0,x)>0$ on $]0,x_M[,$  it implies that $b(a,x,y) (N-N_1)(a,y)$ is of constant sign. We integrate in $x$ and find that $B(a,y) (N-N_1) (a,y)$ is of constant sign. Using assumption (\ref{as:Bpositiv}), $N-N_1$ is of constant sign on $[A^*,A_1]\times ]0,x_M[.$ Formula of characteristics (\ref{eq:char:N}) then gives that $(N-N_1)(0,]Y(A^*,0),Y(A^*,x_M)[)$ is of constant sign. Since $Y(A^*,0)\leqslant 0$ and $Y(A^*,x_M)\geqslant x_M,$ we deduce that it is of constant sign on $\{0\}\times ]0,x_M[,$ thus, using once more the formula of characteristics, on $\R_+^2.$  Since $\iint (N-N_1) dx\, da =0,$ we deduce $N=N_1.$


				\section{Various Generalizations of the Results}


		\subsection{Division Rate With Compact Support}
\label{part:compact}

We have used previously strong assumptions on the support of $b$ and $B.$ Here we relax it and suppose:
\beq \label{as:suppbBcomp}
\exists \;A>0,\quad Supp_a (b,\, B) \subset [0,A].
\eeq
To illustrate this on a simple example, let us take, for $B>0$ constant,
$B(a,x)=B\1_{a\leqslant A}.$
We integrate formula (\ref{def:opG}) and find, for any $f\in E$ such that $\int\limits_0^{x_M} f(x)dx=1:$
$$\int\limits_0^{x_M} \opG_\lb (f)(x) dx= 2 B \int\limits_0^A\int\limits_0^{x_M} 
f(x) e^{-\lb a - Ba} dx\, da = \f{2B}{\lb + B} (1-e^{-(\lb + B)A})=\mu_\lb.$$
It implies that there exists $\lb>0$ 
such that $\mu_\lb = 1$ \emph{iff} $\mu_{\lb=0} >1,$ that is \emph{iff}
$B>\f{1}{A} Log(2).$

As seen on this example, we need a condition to ensure that  $\lb_0 \geqslant 0.$ Namely:

\beq\label{as:comp:H2}
\alpha=\sup\limits_{a,y} \f{\int b(a,x,y) e^{-\int\limits_0^a B\big(s,X(s,x)\big) ds} dx}{ B(a,y)} < \f{1}{2}.
\eeq
(We can divide by $B(a,y)$ without loss of generality since $B(a,y)=0$ implies $b(a,x,y)=0$.)

Such a condition is not surprising and occurs always in structured population equations (see \cite{BP1}). We check easily that it generalizes the condition 
of the example $B(a,x)=B\1_{a\leqslant A}.$

To have uniqueness of a solution and of a $\lb>0,$ we need to replace assumption (\ref{as:bpositiv}). Indeed, we have noticed in part \ref{sec:main} that it generally implies $Supp_a(b)$ unbounded. If we suppose $X(A,x_M)<x_M,$ for $X(A,x_M)\leqslant x < x_M,$ if $y<\min\limits_{0\leqslant a\leqslant A} X(a,x_M),$ 
then for all $0\leqslant a\leqslant A,$ one has $X(a,y)< X^2(A,x_M)\leqslant X(A,x_M) <x,$ so (\ref{as:bpositiv}) cannot be verified. We make the following  restrictive assumption:
\beq\label{as:Gammapositiv}
\forall \; 0<a<x_M,\quad \forall\; 0<x<x_M,\quad \Gamma(a,x)> \,0.
\eeq
Without assumption (\ref{as:Gammapositiv}), we are unable to prove uniqueness of a solution. Refer to the proof of theorem \ref{th:uniqueness}: the fact that $\iint N\phi \,dx\,da >0$ plays a central role, and here we would only be able to prove that $N>0$ on an interval $]0,x_L[$ and $\phi>0$ on $]x_L,x_M[$ where $x_L$ is such that $X(A,x_L)=x_L$ (this proof can be done recursively on the basis of formulae (\ref{form:N}) and (\ref{form:phi})). 
We complete assumption (\ref{as:Bpositiv}) by an assumption on $Supp_a (b):$ 
\beq\label{as:bpositiv:comp}
\forall\; x \in ]0,x_M[,\;\forall \;y\in [x,x_M[,\quad b(A^*\leqslant a\leqslant A_1,x,y)da >0.
\eeq

\begin{theorem}\label{th:comp}
Under assumptions (\ref{as:suppb}), (\ref{as:gamma3})--(\ref{as:b4}), 
(\ref{as:suppbBcomp}), (\ref{as:comp:H2}), there exists solutions $(N,\phi)\in \cac^1 (\R_+^2)$ respectively to problems (\ref{eq:cellcyclin})  (\ref{eq:cellcyclind}) for a $\lb_0 > 0;$ and $Supp_a \phi \subset [0,A].$ 
Under the additional assumptions  (\ref{as:Bpositiv}), (\ref{as:Gammapositiv}) and (\ref{as:bpositiv:comp}), the solution $(\lb_0,\;N,\;\phi)$ of the eigenproblem (\ref{eq:cellcyclin})(\ref{eq:cellcyclind}) is unique and $\lb_0 >0.$
\end{theorem}
\begin{proof}
We proceed in the same way than in part \ref{part:strong}, so we let the details of the proof to the reader. 
We first take $\ep>0$ and consider the following regularization of the parameters, which is the symmetric of definition (\ref{def:bBep:strong}): $b_\ep(a,x,y)=b(a,x,y) + \f{\ep}{x_M}\1_{a\leqslant A},\quad B_\ep (a,y) = B(a,y) + \ep \1_{a\leqslant A}.
$
With this definition, we carry out the same calculations and define the integral operator $\opG_\lb^\ep$ by the same formula (\ref{def:opG}). 
\begin{lemma}\label{lem:comp2}
Under assumptions (\ref{as:suppb}), (\ref{as:gamma3})--(\ref{as:b4}), (\ref{as:suppbBcomp}), for all $\ep \geqslant 0,$ $\lb \in \R,$ the operator $\opG^{\ep}_\lb$ is compact.
\end{lemma}
The proof is the same than for lemma \ref{lem:comp}, even simpler since we do not have any difficulty with the convergence in the $a-$variable of the integrals. For this reason, we can now let $\lb <0$.
\begin{lemma}\label{lem:reg2}
For all $\ep >0$ and $\lb \in \R,$ under assumptions of lemma \ref{lem:comp2}, there exists a unique $\mu_{\lb,\ep} >0$ and a unique $N^0_{\lb,\ep}\in E,$ $N^0_{\lb,\ep} >0,$ such that 
		$\opG_{\lb}^\ep(N^0_{\lb,\ep})=\mu_{\lb,\ep} N^0_{\lb,\ep},\quad
||N^0_{\lb,\ep}||_E = 1.
$
What is more, $\lb \to \mu_{\lb,\ep}$ is a continuous decreasing function which vanishes when $\lb\to +\infty.$ For $\lb=0,$ under assumption (\ref{as:comp:H2}), we have $\mu_{0,\ep} >1$ for $\ep$ small enough. 
\end{lemma}
\begin{proof}
The proof is the same than for lemma \ref{lem:reg}, except that we can let $\lb < 0.$ It only remains to prove that $\mu_{0,\ep}>1.$ 
We denote by $N^0_{0,\ep}\in E$ the (unique up to a multiplicative constant) solution of
$\opG_0^\ep (N^0_{0,\ep})=\mu_{0,\ep} N^0_{0,\ep}.$
Defining $N_{0,\ep}$ by the characteristic formula (\ref{eq:char:N}) with $N_{0,\ep}^0,$ $B_\ep$ and $b_\ep$ instead of $N,$ $B$ and $b,$ the function $N_{0,\ep}$ is then solution to the following problem:
\begin{equation} \left\{
	\begin{array}{cl}
& \f{\p}{\p a} N_{0,\ep} + \f{\p}{\p x}[ \Gamma(a,x) N_{0,\ep} ] + B_\ep(a,x) \, N_{0,\ep} =0, \quad \, a \geqslant 0, \; x \geqslant 0, 
\\ \\
&\mu_{0,\ep} N_{0,\ep}(a=0,x) = 2 \int b_\ep(a,x, y) N_{0,\ep}(a,y) dy \, da,
\\ \\
& N_{0,\ep} \geqslant 0, \qquad \int N_{0,\ep} (a=0,x)dx =1.
\end{array} \right.
\label{eq:cellcyclin:ep:lb}
\end{equation}
Notice that we need $\ep >0$ in order to have convergence of the integral $\int N_{0,\ep} dadx:$ indeed, $Supp_a N_{0,\ep}$ is not compact.
We integrate this equation in $a$ between $0$ and $A$ 
and in $x,$ and find:
$$\int N_{0,\ep} (A,x) dx + (1-\f{2}{\mu_{0,\ep}})\int B_\ep (a,x) N_{0,\ep} (a,x) dx\, da =0.$$
Since $\int B_\ep N_{0,\ep} dx\, da = \f{\mu_{0,\ep}}{2},$ it can be written
$\int N_{0,\ep} (A,x) dx = 1 - \f{\mu_{0,\ep}}{2}.$ 

\noindent 
Hence, $\mu_{0,\ep} > 1$ \emph{iff} $\int N_{0,\ep} (A,x) dx < \f{1}{2}.$
The characteristic formula allows us to write:
$$
\int N_{0,\ep} (A,x) dx = 
\int N_{0,\ep}(0,Y(a,x))e^{-\int\limits_0^A \big\{  \f{\p}{\p x} \Gamma \big(s,X\big(s,Y(a,x)\big)\big) + B_\ep \big(s,X\big(s,Y(a,x)\big)\big) \big\}ds}dx$$
$$
=\int N_{0,\ep}(0,x)e^{ -\int\limits_0^A  B_\ep \big(s,X(s,x)\big)ds} dx
=\f{2}{\mu_{0,\ep}} \int b_\ep(a,x,y) N_{0,\ep} (a,y) e^{ -\int\limits_0^A  B_\ep \big(s,X(s,x)\big)ds}dx\, dy\, da
.
$$
Under assumption (\ref{as:comp:H2}), which is also verified by $b_\ep$ and $B_\ep$ for $\ep$ small enough (to be more precise, if $e^{-\ep A} \leqslant \f{1}{2}$), it implies $\int N_{0,\ep} (A,x)dx \leqslant \f{2}{\mu_{0,\ep}}\f{1}{2}\f{\mu_{0,\ep}}{2}=\f{1}{2},$ and it ends the proof of lemma \ref{lem:reg2}.
\end{proof}

We are now ready to prove theorem \ref{th:comp}. Lemma \ref{lem:reg2} gives us a unique solution $N_\ep^0$ and a unique $\lb_\ep >0$ for which $\mu_{\lb_\ep,\ep}=1$ and $||N_\ep^0||_E = 1.$ 

The family $(N_\ep^0,\lb_\ep)_{0<\ep<1}$ is compact, so we can extract a subsequence tending to a solution $(N^0,\lb_0)$ of the equation
$$N^0(x)=2\int b\big(a,x,X(a,y)\big) N^0(y)e^{-\lb_0 a -\int\limits_0^a B\big(s,X(s,y)\big)ds}\,dy\,da.$$

Defining $N$ by (\ref{eq:char:N}), it remains to check that $\lb_0>0,$ which will imply $\int N dx\, da <\infty.$ 
Since $N^0\geqslant 0$ and $N^0\neq 0,$ we have $\int N^0(x) dx=2\iint B(a,x) N(a,x)\,dx\,da>0,$
so integrating equation (\ref{eq:cellcyclin}) in $x,$ and in $a$ from $0$ to $A,$ one gets:
$$\int N(A,x)dx+\lb_0\int\limits_0^A\int\limits_0^{x_M} N(a,x) dx\, da=\iint B(a,x) N(a,x) dx\, da.$$
Using the preceding calculation done for $\int N_{0,\ep} (A,x) dx$ and assumption (\ref{as:comp:H2}), we deduce:
$$\lb_0\int\limits_0^A\int\limits_0^{x_M} N(a,x) dx\, da=\int \biggl(B(a,y)
- 2\int b(a,x,y)e^{-\lb_0 A-\int\limits_0^A B\big(s,X(s,x)\big)ds}dx\biggr)N(a,y) dy\, da
>0,
$$
which implies $\lb_0 >0.$
To prove uniqueness, 
and to find a solution for problem (\ref{eq:cellcyclind}), the proof is identical to the one of part \ref{sec:proofs}, so we need to prove $\int N \phi \,dx\, da >0.$ 

Since $N\neq 0,$ there exists $y_1\in ]0,x_M[$ such that $N(0,y_1)>0.$ Formula (\ref{eq:char:N}) and assumption (\ref{as:bpositiv:comp}) implies then $N(0,x)>0$ for $x\leqslant X (A_1,y).$ Recursively, it stands for $x\leqslant X^n (A_1,y)$ for all $n\in \mathbb{N}.$ Under assumption (\ref{as:Gammapositiv}), $X$ is increasing, so the sequence $X^n (A_1,y)$ tends to a limit $l$ that verifies $X(A_1,l)=l,$ so $l=x_M$ and $N(0<x<x_M)>0.$ Hence, $\int N\phi\, dx\, da >0.$ 
\end{proof}
 


				\subsection{Weak Theory When $b(a,x,y)$ is Continuous in the $x-$Variable}
\label{part:gen:weak}

In this section, we extend the previous results to a larger class of parameters.
First, we can relax assumption (\ref{as:H1}), and replace it by the two following conditions.
 \beq\label{as:H1bis2}
C_0(x)=\int\int \underbrace{b(a,x,X(a,y))}_{birth \;\;term}\underbrace{e^{-\int_0^a B(s,X(s,y))ds}}_{death \;\;term}dady \in \cac^0 \bigl([0,x_M]\bigr). 
\eeq 
This assumption is necessary to obtain compactness, independently of $\ep\geqslant 0,$ and to apply Ascoli theorem in $E=\cac^0 \bigl([0,x_M]\bigr).$ 
\beq\label{as:H1bis1}
\int_0^\infty B(s,X(s,y))ds>\ln(2),\quad \forall y\in ]0,x_M[.
\eeq 
This condition means that there is enough birth along the characteristic curves. We check as in paragraph \ref{sec:main} that it is verified by the given examples of part \ref{sec:examples}. It is used to prove, with the preceding notations, that the eigenvalue $\mu_{0,\ep} >1$ (as in paragraph \ref{part:compact} it is no more necessarily equal to $2$) and hence that the eigenvalue $\lb_0 >0.$

We also relax the regularity assumptions on the coefficients, and we suppose: 
\beq\label{as:regBGamma:weak}
B \in L^\infty_{loc} (\R_+^2),\qquad \Gamma \in W^{1,\infty} (\R_+^2).
\eeq
We obtain the following theorems.
\begin{theorem}\label{th:cellcyclin:gen}
Under assumptions (\ref{as:suppb}), (\ref{as:H1bis2})--(\ref{as:regBGamma:weak}), with the positivity assumptions (\ref{as:Bpositiv}) and (\ref{as:bpositiv}), there exists unique solutions $(\lb_0,\;N,\;\phi)$ to the eigenproblems (\ref{eq:cellcyclin}) (\ref{eq:cellcyclind}), and $\lb_0 >0.$ Moreover, for all $\;0<\eta<1,$ $Ne^{\lb_0\eta a}$ and $\phi\, e^{\lb_0 \eta a} \in L^1(\R_+\times[0,x_M]).$ 
\end{theorem}

\begin{theorem}\label{th:cellcyclin:gen:comp}
Under assumptions (\ref{as:suppb}), (\ref{as:suppbBcomp}) and (\ref{as:comp:H2}), (\ref{as:H1bis2})--(\ref{as:regBGamma:weak}), with the positivity assumptions (\ref{as:Gammapositiv}) and (\ref{as:bpositiv:comp}), there exists unique solutions $(\lb_0,\;N,\;\phi)$ to the eigenproblems (\ref{eq:cellcyclin}) (\ref{eq:cellcyclind}), and $\lb_0 >0.$ Moreover, for all $0<\eta<1,$ $Ne^{\lb_0\eta a}\in L^1(\R_+\times[0,x_M])$ and $\phi \in L^\infty (\R_+^2),$ $Supp_a \phi \subset [0,A]\times [0,x_M]).$ 
\end{theorem}

\begin{proof}
Refer to proofs of theorems \ref{th:cellcyclin}, \ref{th:cellcyclind}, \ref{th:uniqueness} and \ref{th:comp}.
\end{proof}


\subsection{Case of Equal Repartition After Division}\label{dirac}


In this section, we consider the case
\beq\label{eq:defb}
b(a,x,y)=\delta_{x=\f{y}{2}} B(a,y),
\eeq
With $B\in \cac(\R_+^2)$ and $\Gamma$ verifying the Cauchy-Lipschitz conditions, $B(a,x\geqslant x_M)=0.$ Our proof can be extended to the cases of unequal mitosis, \emph{i.e.} if $b$ is a sum of Dirac measures.
We cannot use lemma \ref{lem:reform} since its formulae are not verified anymore, but the method remains.
To make the details easier to understand, we also make the following assumptions (which is reasonable if we follow \cite{CBBP1}):
\begin{equation}\label{as:xMx0}
\begin{array}{cl}
& \forall \;\; 0<x \leqslant \; \f{x_M}{2},\quad \forall\;a\,\geqslant\,0,\quad \Gamma(a,x)> 0,
\\ \\
& \big\{(a,x)\in \R_+\times ]0,x_M[;\;\quad s.t.\quad \Gamma(a,x)=0\big\}:=x_0(a)\;\text{is an increasing curve}.
\end{array}
\end{equation}
We make assumption (\ref{as:gamma1}) with $\wt \Gamma$ regular. It allows us to divide $\Gamma$ by $x$ and obtain a smooth function.
We first establish the following lemma, which replaces lemma \ref{lem:reform} and on which is based the proof.

\begin{lemma}\label{lem:reform:dirac}
Let $B\in \cac_b (\R_+^2),$ $b$ defined by (\ref{eq:defb}) and $\Gamma\in\cac^1 (\R_+^2)$ verifying assumption (\ref{as:xMx0}). If $N$ is solution of problem (\ref{eq:cellcyclin})(\ref{as:suppN}) then
the two following identities stand:
\beq \label{eq:charN:dirac}
N(0,x)= 4 \int\limits_{\alpha(x)}^\infty B(a,2x) N\big(0,Y(a,2x)\big)e^{-\lb_0 a- \int\limits_0^a (\p_x \Gamma + B)\big(s,X(s,Y(a,2x))\big)ds}da,
\eeq 
where $\alpha(x)$ is defined by $\alpha(x)=Inf\big\{a\geqslant 0,\quad Y(a,2x)\leqslant \f{x_M}{2}\big\}.$ 
\beq \label{eq:charN:int}
N(0,x>0) =  4 \int\limits_{0}^{Min(\f{x_M}{2},2x)} \f{B(f(a,2x),2x)}{\Gamma(f(a,2x),2x)} N(0,a) e^{-\lb_0 f(a,2x) - \int\limits_0^{f(a,2x)} B(s,X(s,a))ds}da\,\1_{x\leqslant \f{x_M}{2}},
\eeq 
where $f(.,2x)$ is the inverse function of $Y(.,2x)$ and is so defined by $ f(Y(a,2x),2x)=a.$
\end{lemma}

\begin{proof}
By the method of characteristics previously used, we obtain by straightforward calculation:
\beq\label{eq:charN:dirac1}
N(0,x)= 4 \int\limits_0^\infty B(a,2x) N\big(0,Y(a,2x)\big)e^{-\lb_0 a- \int\limits_0^a (\p_x \Gamma + B)\big(s,X(s,Y(a,2x))\big)ds}da.
\eeq 
The flow $Y(a,x)$ can be defined also by the following differential system equivalent to (\ref{eq:char}):
\beq \left\{
\begin{array}{cl} &
\f{d}{d a} Y(a,x)= -\Gamma (a, x)e^{-\int\limits_0^a \partial_x \Gamma \big(s,X(s,Y(a,x))\big) ds},\quad a\geqslant 0,\; x\geqslant 0, 
\\ \\
& Y(0,x)=x,\quad x\geqslant 0.
\end{array}
\right. 
\label{eq:char:Y}
\eeq
For $x\geqslant \f{x_M}{2},$  it implies that $\f{d}{da} Y(a,2x) \geqslant 0$ so $Y(a,2x)\geqslant 2x\geqslant x_M.$ Since we are looking for solutions verifying (\ref{as:suppN}), formula (\ref{eq:charN:dirac}) implies that 
$N(0,x\geqslant \f{x_M}{2})=0.$
Hence, in formula (\ref{eq:charN:dirac1}) the integral can be reduced to $a$ such that $Y(a,2x)\leqslant \f{x_M}{2},$ and we find formula (\ref{eq:charN:dirac}).

\

\noindent
For $x\leqslant \f{x_M}{4},$ assumption (\ref{as:xMx0}) implies that $Y(a,2x)$ is decreasing with the $a$ variable, so $Y(a,2x)\leqslant \f{x_M}{2},$  $\alpha(x)=0.$ We make the change of variables
$a\rightarrow a'=Y(a,2x).$
When $a=0,$ $a'=2x$ and when $a=+\infty,$ $a'$ tends to a limit $l\leqslant 0,$ since for $l>0$ one has $\Gamma(a,l)>0$ so $\f{d}{da} Y(a,l) <0.$ Formula (\ref{eq:charN:int}) comes, since we have
$da=\f{d}{da} Y(a,2x) da'=-\Gamma (a, 2x)e^{-\int\limits_0^a \partial_x \Gamma \big(s,X(s,Y(a,2x))\big) ds}da'.
$

\

\noindent
For $\f{x_M}{4}< x \leqslant \f{x_M}{2},$ since $2x > \f{x_M}{2},$ for $a$ small we have $Y(a,2x)\geqslant \f{x_M}{2},$ so $\alpha(x)>0,$ and the definition of $\alpha(x)$ implies $\Gamma(\alpha(x),2x)\geqslant 0.$ 
Because we made assumption (\ref{as:xMx0}), $\Gamma(a> \alpha(x),2x)>0$ so $Y(a>\alpha(x),2x)$ is decreasing and we can define once more the change of variables 
$a\rightarrow a'=Y(a,2x).$
When $a=\alpha(x),$ $a'=\f{x_M}{2}$ and when $a=+\infty,$ $a'$ tends to a limit $l\leqslant 0.$ Formula (\ref{eq:charN:int}) then comes.
\end{proof}

\

We now formulate an equivalent condition of assumption (\ref{as:H1}), used to obtain compactness:
\beq\label{as:H1:dirac}
h(x):=\int\limits_{\alpha(x)}^\infty e^{-\int\limits_0^a \f{\p}{\p x} (\Gamma + B)\big(s,X(s,Y(a,x))\big)ds} \1_{Y(a,x)\in [0,\f{x_M}{2}]} da \in \cac^b_x ([0,x_M]).
\eeq 
Formulated in terms of partial differential equations, it means that the solution $v$ to (\ref{as:H1:PDE}) satisfies $v\in  \cac^b_x\big([0,x_M],L^1_a(\R_+)\big).$
For the adjoint problem, we replace assumption (\ref{as:suppb}) by:
\beq\label{as:suppB:dirac}
B\big(a,x\geqslant X(a,x_M)\big)=0.
\eeq

To have uniqueness, we make the following assumption:
\beq\label{as:Bpositiv:dirac}
\forall\;x\in]0,x_M[,\quad \forall\;a\in\R_+,\quad B(a,x)>0.
\eeq 

\begin{theorem}\label{th:dirac}
 Under the assumptions of lemma \ref{lem:reform:dirac} and the supplementary assumption (\ref{as:H1:dirac}), there exists a solution $(N,\lb)$ to problem (\ref{eq:cellcyclin})(\ref{as:suppN}) and $\lb_0 > 0.$ Under the supplementary assumption (\ref{as:suppB:dirac}), which implies (\ref{as:suppN}), there is a solution $\phi$ to the adjoint problem (\ref{eq:cellcyclind}) related to an eigenvalue $\lb_0>0.$ Under the positivity assumption (\ref{as:Bpositiv:dirac}), $(\lb_0,N,\phi)$ are unique.
\end{theorem}

\begin{proof}
We follow the same way than for the proof of theorem \ref{th:cellcyclin} for instance, so we detail only the specific points.
The main new difficulty is that on formula (\ref{eq:charN:int}), we have divided by $\Gamma\big(f(a,2x),2x\big)$ on a domain where it is strictly positive, so its inverse is well-defined, but it vanishes if $x$ tends to $0$ or $\f{x_M}{2}.$ It does not matter however, since when it vanishes, its inverse is multiplied by $N(0,a)$ which also vanishes. We are then led to use either of formulae (\ref{eq:charN:dirac}) or (\ref{eq:charN:int}) according to where $x$ stands.
 
\

We first define a regularised operator $\opG_{\ep,\lb}:\;X\rightarrow \;X$ on 
$ X=\cac ([0,\f{x_M}{2}])$ 
 by either of the equivalent formulae:
$$\opG_{\ep,\lb}(g)(x)=4 \int\limits_{\alpha(x)}^\infty B(a,2x) g\big(Y(a,2x)\big)e^{-\lb - \int\limits_0^a (\p_x \Gamma + B)\big(s,X(s,Y(a,2x))\big)ds}da+\ep\int\limits_0^{\f{x_M}{2}} g(x)dx$$
$$=4 \int\limits_{0}^{Min(\f{x_M}{2},2x)} \f{B(f(a,2x),2x)}{\Gamma(f(a,2x),2x)} g(a) e^{-\lb f(a,2x) - \int\limits_0^{f(a,2x)} B(s,X(s,a))ds}da+\ep\int\limits_0^{\f{x_M}{2}} g(x)dx.$$ 
The difference with the previous regularisations of parts \ref{part:strong} and \ref{part:compact} stands in the fact that we regularize differently: changing $B$ in $B_\ep$ would not change the value of $\opG(g)(x=0)$ which would remain zero.
\begin{lemma}\label{lem:comp:dirac}
Under the assumptions of lemma \ref{lem:reform:dirac} and the supplementary assumption (\ref{as:H1:dirac}), for all $\ep\geqslant 0,$ $\lb \geqslant 0,$ the operator $\opG_{\ep,\lb}:\; X\rightarrow X$ is compact.
\end{lemma}
\begin{proof}
It suffices to take $\ep=\lb=0.$ 
The assumptions ensure that the operator $\opG=\opG_{0,0}$ is well-defined. To apply Ascoli theorem, let $\eta>0$ arbitrarily small, we look for $\nu>0$ such that 
$$\forall \; g\in X, \;||g||_\infty \leqslant 1,\quad \forall \; x_1,\;x_2\in [0,\f{x_M}{2}],\quad 0< x_1-x_2|< \nu \Rightarrow \Delta (x_1,x_2)=|\opG(g)(x_1)-\opG(g)(x_2)|< \eta.$$ We distinguish three cases: around $\f{x_M}{2},$ around $0,$ and on the compact subset $[\delta,\f{x_M}{2}-\delta]$ with $\delta>0$ small enough.

\begin{enumerate}
 \item 
For $x_1,\;x_2$ close to $\f{x_M}{2},$ we use formula (\ref{eq:charN:dirac}) to define $\opG,$ and remark that $\lim_{x\to \f{x_M}{2}}\alpha(x)=\infty.$
Indeed, for $x$ close to $\f{x_M}{2},$ $\Gamma(0,x)<0$ and under assumption (\ref{as:xMx0}) the function $\Gamma(a,x)$ remains negative till $a=a_0(2x),$ $a_0(2x)$ defined by $\Gamma(a_0(2x),2x)=0.$ 
The curve
 $Y(a,2x)$ defined by (\ref{eq:char:Y}) increases for $a \leqslant a_0(2x),$ so  for $a\leqslant a_0(2x)$ one has $Y(a,2x) \geqslant 2x >\f{x_M}{2}:$ hence, $\alpha(x)>a_0(2x).$
Under assumption (\ref{as:xMx0}) it is clear that $\lim_{2x\to x_M} a_0(2x)=\infty,$ so $\lim_{2x\to x_M} \alpha(x)=\infty.$

Assumption (\ref{as:H1:dirac}) then implies $\lim_{2x\to x_M} \opG(f)(x)= 0$ uniformally in $g,$ if $||g||_\infty \leqslant 1.$ 

\item
For $x_1,$ $x_2>0$ close to $0,$ we have $\opG(g)(0)=0$ since $B(a,0)=0.$ Since $B$ is uniformally continuous, and the integral operator is uniformally convergent thanks to assumption (\ref{as:H1:dirac}), we can bound $\opG(g)(x_{1,2})$ uniformally for $x_{1,2}$ small enough and $||g||_\infty \leqslant 1.$ 

\item 

For ${\delta} < x_1,\;x_2 <\f{x_M}{2}-{\delta},$ we use formula (\ref{eq:charN:int}) to define $\opG.$ This formula gives a classical form of $\opG$ as an integral operator: we know that it is compact as soon as the kernel under the integral is continuous and bounded. It remains to prove that $\Gamma(f(\sigma,2x),2x)$ does not vanish.

For $x<\f{x_0(0)}{2},$ one has $2x <x_0(0)\leqslant x_0(a)$ for all $a\geqslant 0,$ so under assumption (\ref{as:xMx0}) we have $\Gamma(a,2x)>0$ for all $a>0.$

For $\f{x_0(0)}{2} \leqslant x < \f{x_M}{2},$ for all $\sigma$ we have proved above that $f(\sigma,2x)\geqslant \alpha(x) >a_0 (2x),$ which implies by definition of $a_0(x)$ that $\Gamma(f(\sigma,2x),2x)>0.$

Since $\Gamma(f(\sigma,2x),2x)>0$ on the compact subset $(x,\sigma)\in [\f{\delta}{2},\f{x_M}{2}-\f{\delta}{2}]\times [0,Min(2x,\f{x_M}{2})],$ and $\Gamma$ is continuous, it reaches its minimum $\Gamma^{inf} >0:$ this ends the proof of lemma \ref{lem:comp:dirac}.

\end{enumerate}

\end{proof}

\begin{lemma}\label{lem:KR:dirac}
 For all $\ep >0$ and $\lb\geqslant 0,$ under the assumptions of lemma  \ref{lem:comp:dirac}, there exists a unique $\mu_{\ep,\lb} >0$ and a unique $N^0_{\ep,\lb}\in X$ such that 
$\opG_{\ep,\lb}(N^0_{\ep,\lb})=\mu_{\ep,\lb} N^0_{\ep,\lb},\quad ||N^0_{\ep,\lb}||_X =1.$ 
Moreover, $\lb\to\mu_{\ep,\lb}$ is a continuous decreasing function, with $\mu_{\ep,\infty}=\ep\f{x_M}{2}$ and $\mu_{\ep,0}=2+\ep\f{x_M}{2}.$
\end{lemma}
We let the reader check the proof (equivalent to that of lemma \ref{lem:reg}).
We define, as soon as $\ep\f{x_M}{2}<1,$ a unique $\lb_\ep >0$ such that $\mu_{\ep,\lb_\ep}=1.$ We denote as before $N_\ep^0$ the associated eigenvector with $||N^0_\ep||_X =1.$ As in lemma \ref{lem:reg}, it comes from the previous study that the family $(N^0_\ep)_{0\leqslant \ep\leqslant 1}$ is compact, so we extract a subsequence tending to $(\lb_0,N^0)\in \R_+\times E.$ 
Under assumption(\ref{as:Bpositiv:dirac}), formula (\ref{eq:charN:int}) implies that if $N^0(x_1)>0$ then $N^0(x\in ]\f{x_1}{2},\f{x_M}{2}[)>0,$ so recursively it implies $N^0(]0,\f{x_M}{2}[)>0.$ 

The resolution of the adjoint problem is made as in theorem \ref{th:cellcyclind}. The solution of the adjoint problem (\ref{eq:cellcyclind}) can be written now as
$$\phi(0,x)=2\int\limits_0^\infty B\big(a,X(a,x)\big)\phi(\f{X(a,x)}{2})e^{-\lb_0 a -\int\limits_0^a B\big(s,X(s,x)\big)ds}\,da.$$

\noindent
Since $\phi \ne 0$ and $\f{X(a,x)}{2}\leqslant \f{x_M}{2},$ there exists $0<x_1<\f{x_M}{2}$ where $\phi(x_1)> 0.$ So $\iint N\phi\,dx\,da >0,$ which implies $\lb_0>0.$ 
The proof of uniqueness of $N$ and $\phi$ is the same than for theorem \ref{th:uniqueness}.
\end{proof}


				\section{Extensions}


%
		\subsection{Resolution of a Model With Multiple Cyclins}
\label{extension_multiple}

As already mentioned, there is a whole variety of proteins and cyclin/CDK complexes which play a role in the cell cycle, and we can also want to structure the model by the $DNA$ content or the size of the cells, etc.  Hence, that would be useful to include in the model the action of several variables.

Let us suppose that we have $n$ variables playing a role in the cell cycle, and denote them by $x_i$ with $1\leqslant i \leqslant n.$
We write ${\g x} = (x_1,...,x_n),\quad {\g \Gamma}=(\Gamma_1,...,\Gamma_n),\quad |{\g x}|=\sqrt{\sum\limits_{i=1}^n x_i^2},$

\noindent
and we define an order on $\R^n$ thanks to the cone $\R_+^n$ by 
$\quad {\g x}\geqslant {\g y} \quad \iff\quad \forall\;1\leqslant\;i\;\leqslant \; n,\quad x_i\geqslant y_i.$

\noindent
The model (\ref{eq:cellcyclin}) may be generalised by:

\begin{equation} \left\{
	\begin{array}{cl}
& \f{\p}{\p a} N(a,{\g x}) + \grad_{\g x} \cdot [ {\Gamma}(a,{\g x}) N (a,{\g x} ] + \big(\lb_0+B(a,{\g x})\big) \, N(a,{\g x}) =0, \quad \, a \geqslant 0, \; {\g x} \geqslant 0, 
\\ \\
&N(a=0,{\g x}) = 2 \int b(a,{\g x}, {\g y}) N(a,{\g y}) da \, d{\g y},
\\ \\
& N \geqslant 0, \qquad \int N da \, d{\g x} =1.
\end{array} \right.
\label{eq:cellcyclin:N}
\end{equation}
In the same way, the adjoint problem becomes:

\begin{equation} \left\{
	\begin{array}{cl}
& -\f{\p}{\p a} \phi  -   {\g \Gamma}(a,{\g x})  \cdot \grad_{\g x} \phi + \big(\lb_0+B(a, {\g x})\big) \, \phi =2 \int \phi(0,{\g y}) b(a,{\g y},{\g x}) d{\g y}, \quad \, a\geqslant 0, \; {\g x} \geqslant 0, 
\\ \\& \phi \geqslant 0, \qquad \int \phi N da \, d{\g x} =1.
\end{array} \right.
\label{eq:cellcyclind:N}
\end{equation}
We use the notation ${\g \wt \Gamma} = (\wt \Gamma_i)_{1\leqslant i\leqslant n}$ and suppose that
${\Gamma_i}(a,{\g x})= x_i  {\wt \Gamma_i} (a,{\g x}),$
with 
\beq \left\{
\beqa
{ {\g\wt\Gamma}} (a,{\g x})& > 0 \quad \text{for }  {\g x}\approx 0,
\\ \\
{ {\g\wt\Gamma}} (a,{\g x}) &\leqslant 0 \quad \text{for }  {\g x} \geqslant  {\g x}_M=(x_{1M},...,x_{nM})
\eeqa 
\right. 
\label{as:gamma1:N}
\eeq
Concerning division, we have the same relations (\ref{as:b1}) and (\ref{as:b2}) than previously seen. Characteristics and their inverse flow are still defined by (\ref{eq:char}) and (\ref{eq:chard}). Assumption (\ref{as:gamma1:N}) implies that for all
$0\leqslant {\g x}\leqslant {\g x}_M,$ $X(a,x)\leqslant\;{\g x}_M,$
so we add the same conditions (\ref{as:suppN}) and (\ref{as:suppb}) to problem (\ref{eq:cellcyclin:N}).
\begin{theorem}
Replacing $x,$ $y\in\R_+$ by $\g x$ and $\g y \in \R_+^n,$  replacing $\f{\p \Gamma}{\p x}$ by $\grad\cdot {\g \Gamma},$ and $y^{\gamma_1}$ by $\prod\limits_{i=1}^n y_i^{\gamma_i}$ with $\gamma_i\geqslant 1,$ all the preceding theorems  \ref{th:cellcyclin}, \ref{th:cellcyclind}, \ref{th:uniqueness}, \ref{th:comp} and \ref{th:dirac} extend to the $n-$dimensional case under the equivalent assumptions in $n$ dimensions. 
\end{theorem}


\subsection{Asymptotic Behaviour of the Linear Evolution Problem}
Having solved the eigenvalue problem, we are now able to characterize the asymptotic behaviour of the solution $n$ to the time-dependent problem (\ref{eq:cellcyclin:t1}).
First, we establish a General Relative Entropy Inequality, using the same formalism than in \cite{MMP2} and \cite{MMP}.
\begin{proposition}
\label{prop:GRE}
Let $p(t,a,x)\geqslant 0,$ $n(t,a,x)$ smooth solutions of problem (\ref{eq:cellcyclin:t1})
and $\Phi(t,a,x) \geqslant 0$ smooth solution of the adjoint problem:
\begin{equation} \left\{
	\begin{array}{cl}
& -\f{\p}{\p t}\Phi - \f{\p}{\p a} \Phi - \Gamma(a,x) \f{\p}{\p x}  \Phi  + B(a,x) \, \Phi =2\int \Phi(t,0,y)\,b(a,y,x)\,dy, \quad \, a \geqslant 0, \; x \geqslant 0, 
\\ \\
& \Phi \geqslant 0.
\end{array} \right.
\label{eq:cellcyclind:t}
\end{equation}
Then we have, for any function $H:$
\begin{equation*}\label{eq:GRE:ext}\begin{array}{cc}
\f{d}{dt}  {\cal H} (\f{n}{p}) (t) = \f{d}{dt} \int p(t,a,x) \Phi (t,a,x) \1_{p(t,a,x)\neq 0} H\big(\f{n(t,a,x)}{p(t,a,x)}\big) da dx=
2 \int p(t,a,x) \Phi(t,0,y) b(a,y,x) 
\\ \\
\biggl[ H\big(\f{n(t,0,y)}{p(t,0,y)}\big)\1_{p(t,0,y)\neq 0} - H \big( \f{n(t,a,x)}{p(t,a,x)}\big)\1_{p(t,a,x)\neq 0}\;+
H'\big(\f{n(t,0,y)}{p(t,0,y)}\big)\1_{p(t,0,y)\neq 0}\bigl(\f{n(t,a,x)}{p(t,a,x)}-\f{n(t,0,y)}{p(t,0,y)}\bigr)\biggr] dadxdy. \\ \\ \end{array}
\end{equation*}

i) If $p>0$ never vanishes, for $H$ convex, we get
$\f{d}{dt}  {\cal H} (\f{n}{p}) (t) \leqslant 0,$

and if $H$ is strictly convex, $\f{d}{dt}  {\cal H} (\f{n}{p}) (t) = 0$ \emph{iff} $n\equiv Cp$ on $Supp(b)$ with $C\geqslant 0$ constant.

ii) If $p\geqslant 0,$ for $H$ convex, positive and non-decreasing, we also get
$\f{d}{dt}  {\cal H} (\f{n}{p}) (t) \leqslant 0.$
\end{proposition}
\begin{proof}
If $p>0$ everywhere, $\1_{p\neq 0} \equiv 1$ and we have by straightforward computation:
$$\f{\p}{\p t}  \bigl(p \Phi H(\f{n}{p})\bigr)  +  \f{\p}{\p a} \bigl(p \Phi H(\f{n}{p})\bigr) + \f{\p}{\p x}[ \Gamma(a,x) p \Phi H(\f{n}{p}) ] = - 2 p H(\f{n}{p}) \int \Phi (t,0,y) b(a,y,x) dy.$$
Integrating this identity and denoting 
${\cal H} (\f{n}{p}) (t) = \int p(t,a,x) \Phi (t,a,x) H\big(\f{n(t,a,x)}{p(t,a,x)}\big) da dx,$
we get:
\begin{multline*}\f{d{\cal H}}{dt} =
 2 \int p(t,a,x) \Phi(t,0,y) b(a,y,x) \big[ H\big(\f{n(t,0,y)}{p(t,0,y)}\big) - H \big( \f{n(t,a,x)}{p(t,a,x)}\big)\;
\\ \\
+ H'\big(\f{n(t,0,y)}{p(t,0,y)}\bigr)\bigl(\f{n(t,a,x)}{p(t,a,x)}-\f{n(t,0,y)}{p(t,0,y)}\bigr)\bigr] dadxdy,
\end{multline*}
because we have noticed that
$$\int H'\big(\f{n(t,0,y)}{p(t,0,y)}\big)\big[\f{n(t,a,x)}{p(t,a,x)}-\f{n(t,0,y)}{p(t,0,y)}\big]p(t,a,x)
\Phi(t,0,y) b(a,y,x)\,da\,dx\,dy=0.$$
If $p\geqslant 0,$ we replace $H(\f{n}{p})$ by $H(\f{n}{p})\1_{p\neq 0}$ in the calculation, since 
$\p_z \big(p\phi H(\f{n}{p})\1_{p\neq 0}\big)=\1_{p\neq 0}\p_z\big(p\phi H(\f{n}{p})\big)$
and 
$$\int H'\big(\f{n(t,0,y)}{p(t,0,y)}\big)\1_{p(t,0,y)\neq 0}\big[\f{n(t,a,x)}{p(t,a,x)}-\f{n(t,0,y)}{p(t,0,y)}\big]p(t,a,x)
\Phi(t,0,y) b(a,y,x)\,da\,dx\,dy=0.$$
\end{proof}

\begin{theorem}\label{th:evol:exist}
Let $n_0 \in L^1(\R_+^;\phi(a,x)dxda),$ $Supp(n_0) \subset \R_+ \times [0,x_M],$ $\Gamma$ satisfying the Cauchy-Lipschitz conditions and $B\in L^\infty (\R_+^2).$

There exists a unique solution $\wt n \in \cac \big(\R_+;\, L^1(\R_+^2)\big)$ of the following problem:
\beq \label{eq:evol}
\left\{
	\begin{array}{cl}
& \f{\p}{\p t}  \wt n +  \f{\p}{\p a} \wt n + \f{\p}{\p x}[ \Gamma(a,x) \wt n ] + (B(a,x) +\lb_0) \, \wt n =0, \quad \, a \geqslant 0, \; x \geqslant 0, 
\\ \\
&\wt n(t,a=0,x) = 2 \int b(a,x, y) \wt n(t,a,y) dy \, da,
\\ \\
& \wt n(t=0, a,x)=n^0 (a,x). 
\end{array} \right.
\eeq 
Moreover, we have the following inequalities, if
$|n^0(a,x)|\leqslant C_0 N(a,x),$

(i) $\forall \;t\geqslant 0,\quad |\wt n (t,a,x)| \leqslant C_0 N(a,x),$

(ii) $n^0_1 \leqslant n^0_2\quad \Rightarrow\quad \wt n_1 (t,x) \leqslant \wt n_2 (t,x),$

(iii) $\int\limits_0^\infty \wt n (t,x) \phi (a,x) dx da = \int\limits_0^\infty n^0 (a,x) \phi(a,x) dx da,$

(iv)  $\int\limits_0^\infty |\wt n (t,x)| \phi (a,x) da dx = \int\limits_0^\infty |n^0 (a,x)| \phi(a,x) dx da.$

\end{theorem}
\begin{theorem}\label{th:evol:reg}
Under the assumptions of theorem \ref{th:evol:exist}, and if we suppose also that $n^0$ satisfies 
$$|n^0 (a,x)| \leqslant N(a,x),\quad |\p_a n^0 (a,x) +\p_x (\Gamma n^0) (a,x)| \leqslant C_1 N(a,x),$$
the solution to (\ref{eq:evol}) satisfies also
$|\p_t \wt n (t,a,x)|\leqslant (C_1 + \lb_0 + ||B||_{L^\infty})N(a,x).$
\end{theorem}
\begin{theorem}\label{th:conv}
Under the assumptions of theorem \ref{th:evol:exist}, and either those of theorem \ref{th:cellcyclin:gen} or of theorem \ref{th:cellcyclin:gen:comp}, defining $m^0=\int n^0(a,x) \phi (a,x) dx\, da,$ the solution to (\ref{eq:evol}) satisfies:
$$\int |\wt n(t,a,x) - m^0 N(a,x)| \phi(a,x) dx\, da \downarrow 0\quad as \quad t\to\infty.$$
\end{theorem}
\begin{proof}
Since the proofs are the adaptation to our model of those of theorems 3.1, 3.2 and 3.4. of \cite{BP1} or of theorems 4.3 and 3.2 of \cite{MMP}, we let the reader check them.
\end{proof}

		\subsection{Application to a Two Phase Model}\label{part:appli:2phase}
As previously seen in part \ref{sec:pres}, equation (\ref{eq:cellcyclin}) can be considered as a simplification of the linearised eigenvalue problem for the two-compartment model described in \cite{CBBP1}:

\begin{equation} 
\label{eq:cellcyclin:2p}
\left\{
	\begin{array}{cl}
& \f{\p}{\p a} P + \f{\p}{\p x}[ \Gamma(a,x) P ] + \big(\lb+B(a,x)+d_1+L(a,x)\big) \, P - \wt G Q =0, \quad \, a \geqslant 0, \; x \geqslant 0, 
\\ \\
& (\lb + \wt G + d_2) Q = L(a,x) P,\quad \, a \geqslant 0, \; x \geqslant 0, 
\\ \\
&P(a=0,x) = 2 \int b(a,x, y) P(a,y) dy \, da,
\\ \\
& P,\;Q \geqslant 0, \qquad \int P\,+\,Q\, dx \, da =1.
\end{array} \right.
\end{equation}
The adjoint problem is:
\begin{equation} 
\label{eq:cellcyclind:2p}
\left\{
	\begin{array}{cl}
&- \f{\p}{\p a} \phi - \Gamma(a,x)\f{\p \phi}{\p x}  + \big(\lb+B(a,x)+d_1+L(a,x)\big) \, \phi - L(a,x) \psi =
 2 \int \phi(0,y) b(a,y, x)\, dy,
\\ \\
& (\lb + \wt G + d_2) \psi = \wt G \phi,\quad \, a \geqslant 0, \; x \geqslant 0, 
\\ \\
& \phi,\;\psi \geqslant 0, \qquad \int \phi P\,+\,\psi Q\, dx \, da =1.
\end{array} \right.
\end{equation}
Here $P$ and $Q$ denote respectively the proliferative and quiescent populations of cells, $\wt G$ is the recruitement function, $L(a,x)$ the number of cells going from the proliferative to the quiescent compartment, and $d_1,$ $d_2$ are the death rates of each population. 
For the sake of simplicity, we have limited our study here to the case when $L$ is constant, and we make the following assumptions for the coefficients.
\beq\label{as:d1d2LG}
d_1>0,\quad d_2\geqslant 0,\quad 
L(a,x)=L\geqslant\, 0,\quad \wt G \geqslant 0.
\eeq
\begin{theorem}\label{th:exist:2phase}
Under the assumptions of theorem \ref{th:uniqueness}, and with the supplementary assumption (\ref{as:d1d2LG}), there exists unique solutions $(P,Q)\in E^2,$ $(\phi,\psi)\in E^2$ to problem (\ref{eq:cellcyclin:2p}) and (\ref{eq:cellcyclind:2p}) for a unique $\lb \in \R.$ Moreover, denoting $\lb_0$ the eigenvalue of problem (\ref{eq:cellcyclin}) and (\ref{eq:cellcyclind}), we have:
\beq\label{eq:linklblb0}
\lb_0=\lb + d_1 + L\f{\lb + d_2}{\lb+\wt G+ d_2}>\,0,
\eeq
or also, defining $G_+=\wt G + d_2,$ $d_+=d_1 - \lb_0$ and $L_+=L+d_1 -\lb_0:$ 
\beq\label{eq:deflb}
\lb=\f{-(G_++L_+)+{\sqrt{(G_++L_+)^2-4(d_+G_++d_2 L)}}}{2}.
\eeq
The following estimate stands for $\lb,$ with $L_i=L+d_1:$
\beq\label{ineq:lb2}
\lb > \underline{\lb}:=\f{1}{2}\biggl(-(G_+ + L_i)+\sqrt{(G_+ + L_i)^2 - 4 (d_1 G_+  + L d_2)}\biggr).
\eeq
\end{theorem}

\begin{proof}
We reduce the system to a single equation on $P,$ since we can write
$Q=\f{L}{\lb + G_+ } P.$
Writing $\lb_0=\lb + d_1 + L\f{\lb + d_2}{\lb+G_+},$
we find:
\beq
\left\{
	\begin{array}{cl}
& \f{\p}{\p a} P + \f{\p}{\p x}[ \Gamma(a,x) P ] + \big(\lb_0+B(a,x)\big) \, P =0, \quad \, a \geqslant 0, \; x \geqslant 0, 
\\ \\
&P(a=0,x) = 2 \int b(a,x, y) P(a,y) dy \, da,
\\ \\
& P \geqslant 0, \qquad \int P\,\big(1+\f{L}{\lb+G_+}\big)\, dx \, da =1.
\end{array} \right.
\eeq
We only have to change the normalization and the previous study applies to this case: theorems \ref{th:cellcyclin} to \ref{th:uniqueness} give us a unique solution $(P,Q,\phi,\psi,\lb_0>0).$
We denote the eigenvalue $\lb_0=f(\lb)$ where $f$ is a continuous increasing function with a unique singularity for $\lb=- G_+.$ 
Since $Q=\f{L}{\lb+G_+} P,$ we need to have $\lb>-\lb_1,$ so for each $\lb_0\in\R,$ there exists a unique convenient $\lb\in]-G_+,+\infty[$ which is given by formula (\ref{eq:deflb}). The term under the square root is always nonnegative since 
$$(G_++L_+)^2-4d_+ G_++d_2 L\big)\geqslant (G_++L_+)^2-4L_+ G_+
=(G_+ -L_+)^2\geqslant 0.$$
The inequality (\ref{ineq:lb}) is given by the fact $\lb > \lb_2$ where $f(\lb_2)=0.$
\end{proof}

\ 

\noindent
{\bf Discussion.} From theorem \ref{th:exist:2phase} we deduce as in theorem \ref{th:conv} the long-time convergence of the solution of the linearised problem (\ref{eq:cellcyclin:t:2p}) towards $(P,Q)e^{\lb t}.$ However, what is experimentally observed is either convergence towards a steady state or exponential growth but only in the early stages (Gompertzian growth: cf. \cite{GW1} and references therein). In \cite{BASGAB} polynomial infinite growth for 15 cell lines is shown, and in \cite{DH} a single-cell model is built, which is able to exhibit such behaviours.
The linearised problem cannot take into account such phenomena, due to  feedback answer or saturation effect:
it can only come from a non-linearity of the model.

But if the linear renewal equation has a relatively simple asymptotic behaviour, the theory for nonlinear models is much more complicated. Several behaviours are possible: chaotic, periodic, convergence towards stable steady states (for recent references on nonlinear population models, see for instance \cite{BBA}, \cite{CC}, \cite{CI}, \cite{CS}, \cite{F1}, \cite{F2}, \cite{M:McKendrick}, \cite{MPR} or \cite{PT}).

In \cite{CBBP1}, following \cite{GW1} and \cite{GW2}, it is proposed that the non linearity comes from the term $G(N(t))$ where the weighted total population $N(t)$ is defined by
$N(t)=\int\limits_0^\infty\int\limits_0^\infty [\phi^*(a,x)p(t,a,x)+\psi^*(a,x)q(t,a,x)]dx\,da,$
where $\phi^*,\;\psi^*$ are given weights. The recruitment function $G$ is taken equal to
\beq\label{eq:defG}
G\big(N(t)\big)=\f{\alpha_1\theta^n+\alpha_2N^n}{\theta^n+N^n},\quad 0<\alpha_2 <\alpha_1.
\eeq
To study the behaviour of the model, the method of \cite{CBBP1} is inspired of the principles of General Relative Entropy. It is based on estimates and on the construction of subsolutions and supersolutions of a quantity $$S(t)=\int\limits_0^\infty\int\limits_0^\infty [\phi(a,x)p(t,a,x)+\psi(a,x)q(t,a,x)]dadx,$$ 
where $(\phi,\psi)$ is the solution of the adjoint linearised eigenproblem (\ref{eq:cellcyclind:2p}) for a proper value of $\wt G=G(N)$  (see also \cite{CCP} for the application of this method to another model).
Proposition 2.7. of \cite{CBBP1} shows that unlimited growth can be obtained under the two following conditions.

\noindent
(H9) For all $N<\infty,$ the eigenvalue $\lb(N)$ corresponding to $\wt{G} = G(N)$ in (\ref{eq:cellcyclin:2p}) satisfies $\lb(N)>0.$

\noindent
(H10) For each corresponding solutions to the systems (\ref{eq:cellcyclind:2p}) with $\wt G= G(N)$, denoted $(\phi_N,\psi_N),$ there exists a uniform constant $C_u$ such that $\phi^*\geqslant C_u\phi_N$ and $\psi^*\geqslant C_u\psi_N.$

\

Proposition 2.5. of \cite{CBBP1} shows that subpolynomial growth can be obtained under the following conditions, if $d_2>0$ and $\alpha_2>0.$ 
 
\noindent
(H7) For $\tilde{G} = G(\infty)=\alpha_2>0,$ the first eigenvalue of (\ref{eq:cellcyclin:2p}) is $\lb(\infty)=0.$

\noindent
(H8) For the corresponding solutions to (\ref{eq:cellcyclin:2p}) and (\ref{eq:cellcyclind:2p}), denoted respectively $(P_2,Q_2)$ and $(\phi_2,\psi_2),$ there exists positive constants $C_2$ and $C_3$ such that $C_3\phi_2\leqslant \phi^*\leqslant C_2\phi_2$ and $C_3\psi_2\leqslant\psi^*\leqslant C_2\psi_2.$

\

But to obtain exactly (H7), it has been observed in the numerical simulations of \cite{CBBP1} that all the parameters $d_1,\,d_2,\,L,\,\alpha_2,\,\alpha_1$ had to be related and chosen very carefully: a very small change in one of the parameters implies that $\lb(\infty)\ne 0,$ and the system either is bounded (if $\lb(\infty)<0$) or grows exponentially (if $\lb(\infty)>0.$)  This can also be seen by
taking a closer look to formula (\ref{eq:deflb}): indeed, we can also write it as
\beq\label{eq:deflb2}
\lb=\f{-2\big(G_+ d_++Ld_2\big)}{ G_++L_++\sqrt{(G_++ L_+)^2 - 4 G_+d_+}}.
\eeq
Since the denominator of this formula is always positive and bounded, one has $\lb=0$ \emph{iff} 
\beq\label{eq:lb=0}
G_+ d_+=-Ld_2.
\eeq 
The simulations carried out in \cite{CBBP1} were all done with $d_2>0,$ so this formula is verified punctually, for special values of the coefficients $\wt G,\,d_2,\,d_1,\,d_2$ and $L$ linked by (\ref{eq:lb=0}).

\

From a biological point of view, this obligation to have coefficients linked by such a relation seems hardly justified. 
But we can also assume $d_2=0:$ from a biological point of view it can be verified for some kinds of cell populations - stem cells for instance: in the quiescent compartment indeed, there is no reason why the cells should die (see \cite{ACR} or \cite{GW1}, and the references therein).

In this case, and supposing also that $\lim\limits_{N\to\infty} G(N)=0$ (which is indeed realistic) we see that condition (\ref{eq:lb=0}) will be always verified: 
we can now obtain a ``robust'' subpolynomial growth - I call ``robust'' a subpolynomial growth which remains true for a whole range of parameters $d_1,\,L,\,n,\,\alpha_1.$ 
This is expressed by the following proposition.

\begin{proposition}\label{prop:pol}
Let us suppose $d_2=0,$ $G(N)$ defined by (\ref{eq:defG}) with $\alpha_2=0,$  $\Gamma(a,x),\,B(a,x)$ given functions verifying the assumptions of one of the theorems \ref{th:uniqueness}, \ref{th:comp} or \ref{th:dirac}.  We denote $\lb_0>0$ the eigenvalue of (\ref{eq:cellcyclin}), and suppose that $0<d_1<\lb_0$ and that $L>\lb_0-d_1.$  The case $\lim\limits_{N\to\infty}G(N)=\wt G=0$ can be represented by the following system:
\begin{equation}
\label{eq:cellcyclin:2p:inf}
\left\{
	\begin{array}{cl}
& \f{\p}{\p a} P_2 + \f{\p}{\p x}[ \Gamma(a,x) P_2 ] + \big(B(a,x)+d_1+L\big) \, P_2 - {\cal Q}_2 =0, \quad \, a \geqslant 0, \; x \geqslant 0, 
\\ \\
& {\cal Q}_2 = (L+d_1-\lb_0)P_2,\quad \, a \geqslant 0, \; x \geqslant 0, 
\\ \\
&P_2(a=0,x) = 2 \int b(a,x, y) P_2(a,y) dy \, da,
\\ \\
& P_2,\;{\cal Q}_2 \geqslant 0, \qquad \int P_2\,+\,{\cal Q}_2\, dx \, da =1.
\end{array} \right.
\end{equation}
The adjoint problem is:
\begin{equation} 
\label{eq:cellcyclind:2p:inf}
\left\{
	\begin{array}{cl}
&- \f{\p}{\p a} \phi_2 - \Gamma(a,x)\f{\p \phi_2}{\p x}  + \big(B(a,x)+d_1+L\big) \, \phi_2 - L \psi_2 =
 2 \int \phi_2(0,y) b(a,y, x)\, dy,
\\ \\
& L \psi_2 = (L+d_1-\lb_0)  \phi_2,\quad \, a \geqslant 0, \; x \geqslant 0, 
\\ \\
& \phi_2,\;\psi_2 \geqslant 0, \qquad \int \phi_2 P_2\,+\,\psi_2 {\cal Q}_2\, dx \, da =1.
\end{array} \right.
\end{equation}
Under assumption (H10),  and under assumption (H8) adapted for the problems(\ref{eq:cellcyclin:2p:inf}) (\ref{eq:cellcyclind:2p:inf}), there exists a constant $C>0$ such that:
$$N(t)\leqslant C t^{\f{1}{n}},\quad \quad \lim\limits_{t\to\infty} N(t)=+\infty.$$

\end{proposition}
\begin{proof}
Since $d_1 <\lb_0$ formula (\ref{eq:deflb2}) with $d_2=0$ implies that (H9) is verified, so proposition 2.7 of \cite{CBBP1} can be applied and proves unlimited growth. 

The proof of the subpolynomial growth is based on the same tools than in \cite{CBBP1}, proposition 2.5: the only difference is that since $d_2=0,$ at infinity we have $\lb(\infty)=G(\infty)=0,$ so the equation for $Q$ in (\ref{eq:cellcyclin:2p}) only expresses that  $Q$  tends to infinity whereas $(\lb+\wt G)Q$ remains finite. We will obtain a relevant problem by an asymptotic analysis: 
formula (\ref{eq:deflb2}) can be written, if $\wt G,\,\lb\to\, 0,$
$\lb\approx \f{\wt G (\lb_0-d_1)}{L+d_1-\lb_0}.$
We can replace it in the second equation of (\ref{eq:cellcyclin:2p}) and find:
$\wt G Q \f{L}{L+d_1-\lb_0} =LP.
$
Noting ${\cal Q}_2=\wt G Q,$ we obtain a problem for the couple $(P_2,{\cal Q}_2)$ which remains meaningfull if $\lb,\wt G$ vanishes, and by choosing an appropriate normalisation its limit is (\ref{eq:cellcyclin:2p:inf}), which adjoint is (\ref{eq:cellcyclind:2p:inf}).

We define
$S_2(t)=\int\limits_0^\infty\int\limits_0^\infty [\phi_2(a,x)p(t,a,x)+\psi_2(a,x)q(t,a,x)]dadx,$
and we calculate
$$\f{dS_2(t)}{dt} (t)= \iint \biggl\{\phi_2[-(L+B+d_1)p + G(N(t))q - \f{\p p}{\p a} - \f{\p}{\p x}(\Gamma_1 p)] + \psi_2(Lp-G(N(t))q)\biggr\} \,dx\, da,$$
$$\f{dS_2(t)}{dt} (t)
= G
\iint \big\{\phi_2(a,x)-\psi_2(a,x)\big\} q(t,a,x) dxda
=G
\iint \f{\lb_0-d_1}{L+d_1-\lb_0}\psi_2 q dxda  
\leqslant C G
S_2(t).$$
Since $N(t)\geqslant C_3 S_2(t)$ thanks to (H8) we obtain
\beq
\label{eq:ineqS2} 
\f{dS_2 (t)}{dt}\leqslant C S_2(t)\f{\alpha_1 \theta^n}{\theta^n+(C_3 S_2(t))^n} 
\eeq
Since it has been proved in \cite{CBBP1} that $\Sigma (t)=a(t+t_0)^{\f{1}{n}},$ for $a$ large enough, is a supersolution to (\ref{eq:ineqS2}), so for $t_0$ such that $\Sigma(0)\geqslant S_2(0),$ we have by the comparison principle
$S_2(t)\leqslant \Sigma(t).$
It ends the proof. 
\end{proof}
 
\begin{figure}[ht]
\begin{minipage}[b]{0.46\linewidth}
\includegraphics[width=\linewidth, height=7cm]{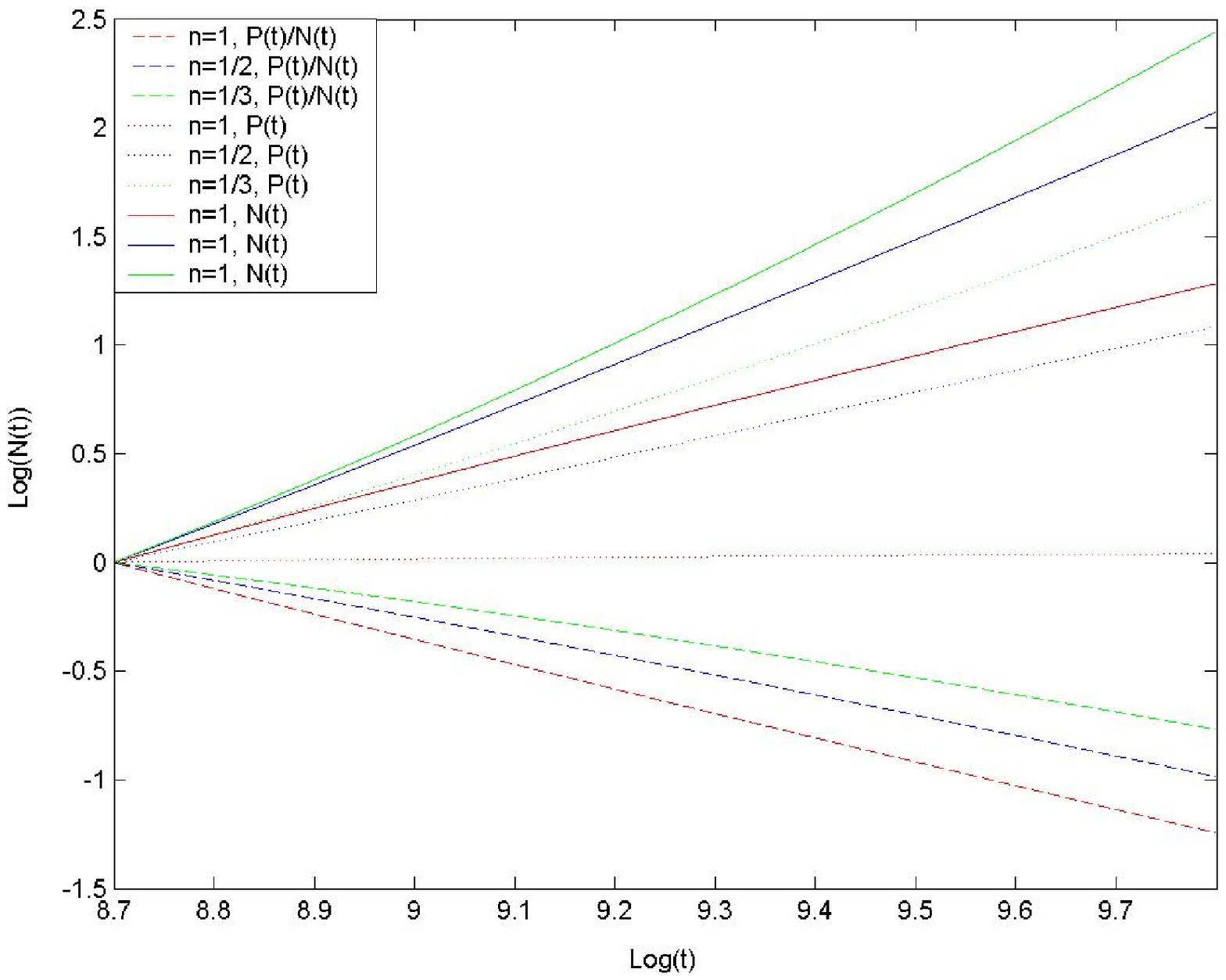} 
\end{minipage} \hfill
\begin{minipage}[b]{0.46\linewidth}
\includegraphics[width=\linewidth, height=7cm]{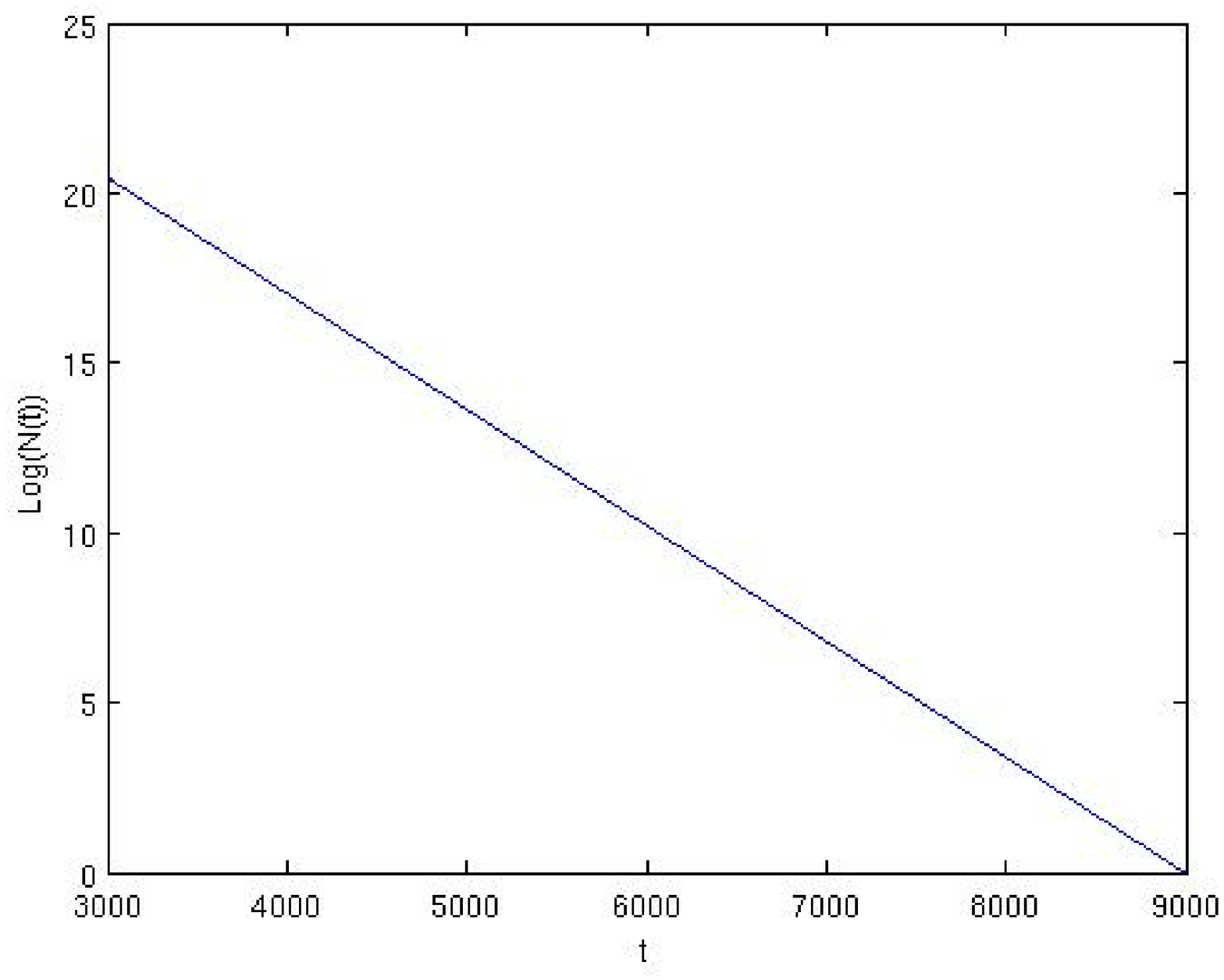}
\end{minipage}\hfill
\caption{\label{fig:sim} Evolution of the total cell population $\iint q(t,a,x) \,dx\,da + \iint p(t,a,x)\,dx\,da$ for a tumoral tissue. Left: with $d_1=0.01,$ with different values of $n=1$ (lower solid line curve, left), $n=1/2$ (medium solid line, left), $n=1/3$ (upper solid line, left) at a Log-Log scale, for $N(t)$ (solid lines) $P(t)$ (dotted lines) and $P(t)/N(t)$ (dashed lines). Right: with $d_1=0.05,$ we see exponential decreasing (Log scale, the three curves for $N(t)$ are superimposed).} 
\end{figure}

For the numerical simulations (see figure \ref{fig:sim}), we take the same values of the parameters than in \cite{CBBP1}, except that $d_2=0:$ we define $\Gamma(a,x)$ by
(\ref{def:gamma}) with
$c_1=0.1,$ $c_2=0.075,$ $r_1=3,$ $c_4=0.4,$ $r_2=1.95.$
We take $\alpha_1=8,$ $\theta=1$ $n=\f{1}{k}$ with $k=1,2,3$ in the definition (\ref{eq:defG}).
We define $L(a,x)$ as in \cite{CBBP1} by:
$$L(a,x)=A_3 \f{A_2^{\gamma_2}}{A_2^{\gamma_2}+x^{\gamma_2}}\1_{[\bar{A},+\infty[}(a),\quad\text{with }\gamma_2=5,\quad A_3=4\quad A_2=2\quad \bar{A}=18.$$
We define $B(a,x)$ by
$B(a,x)=\f{k_1y^{\gamma_1}}{k_2^{\gamma_1}+y^{\gamma_1}}\1_{[A^*,\infty[}(a)$
with
$k_1=1.2\quad k_2=1.5 \quad\gamma_1=5\quad A^*=23.$
We obtain polynomial growth by taking for instance $d_1=0.01$ and it remains true if we make small changes of any coefficient. But if $d_1$ becomes too big, for instance if  we take $d_1=0.05,$ we obtain solutions exponentially vanishing. This is explained by proposition \ref{prop:pol}: indeed, when $d_1$ increases, it becomes bigger than $\lb_0$ and formula (\ref{eq:deflb2}) shows that it implies $\lb(N)<0$ for all $N.$ 

\

{\bf Conclusion.} 
In this article, we first solved the eigenvalue problem (\ref{eq:cellcyclin}) under fairly general assumptions, provided the continuity of the repartition function $b$ or, in the case of equal repartition after division, of the birth rate $B$ (the generalisation to $L^2$ coefficients is a work in progress). Using General Relative Entropy Inequality, we deduced from it the asymptotic behaviour of problem (\ref{eq:cellcyclin:t1}). 

We then applied these results, in part \ref{part:appli:2phase}, to the study of a non linear two cell-compartment model given by equation (\ref{eq:cellcyclin:t:2p}), model which was first introduced by F. Bekkal Brikci, J. Clairambault and B. Perthame in \cite{CBBP1} and \cite{CBBP2} to study the action of proteins on the cell cycle (it can also be considered as a generalisation of the pure size-structured two-compartment models studied by M. Gyllenberg and G.F. Webb in \cite{GW1} and \cite{GW2}.) 

Finally, we exhibited a case of ``robust'' polynomial growth, which reveals coherent with the results of \cite{BASGAB} and \cite{DH}. 

\

This last result could lead to two biological interpretations.
First, when the population becomes larger, the formulation of the eigenvalue problem (\ref{eq:cellcyclin:2p:inf}) (\ref{eq:cellcyclind:2p:inf}), where $Q$ had to be replaced by ${\cal Q}_2$ the limit of $G(N) Q,$ seems to show that the number of quiescent cells tends to infinity more rapidly than the number of proliferating cells: indeed, it seems that $P(t)\approx G(N(t)) Q(t) \approx t^{\f{1}{n}-n},$ so the relative number of proliferating cells, given by $R(t)=\f{P(t)}{N(t)},$ seems to vanish like $t^{-n}.$ It is coherent with the results of \cite{GW1} and it is also confirmed by numerical tests (see figure \ref{fig:sim}, left side). From a biological point of view, it seems true in many cases, for instance for stem cells (see \cite{ACR}), that most of the cells are quiescent. 
Second, this qualitative study seems to emphasize the crucial importance of  apoptosis for homeostasis or tumour growth (see \cite{apoptosis}): indeed, if $d_1=\lb_0$ is reached, the number of tumour cells would decrease rapidly instead of growing to infinity. 

\

{\bf Acknowledgment.} The author expresses very grateful thanks to Beno\^it Perthame for his precious help, ideas and corrections, to Jean Clairambault for many discussions and help in formulating the main biological issues, and to Philippe Michel for his corrections, and for improvements of section \ref{part:gen:weak}.

%
%
%

\end{document}